\documentclass[11pt]{amsart}
\usepackage{amsmath,amsthm,amssymb,appendix,bm,doi,graphicx,hyperref,mathrsfs,geometry,graphics,lipsum,tikz,fancyhdr,etoolbox,lmodern,fontenc,picinpar,verbatim}
\usepackage[capitalise]{cleveref}
\usepackage[OT2,OT1]{fontenc}
\usepackage[all]{xy}
\usepackage[framemethod=tikz]{mdframed}

\usepackage{tikz-cd}

\usepackage{indentfirst}

\newcommand{\mA}{\mathbb A}
\newcommand{\mC}{\mathbb C}
\newcommand{\mP}{\mathbb P}
\newcommand{\mQ}{\mathbb Q}
\newcommand{\mZ}{\mathbb Z}
\newcommand{\cA}{\mathcal A}
\newcommand{\cE}{\mathcal E}
\newcommand{\cM}{\mathcal M}
\newcommand{\cO}{\mathcal O}
\newcommand{\cT}{\mathcal T}
\newcommand{\cV}{\mathcal V}
\newcommand{\cW}{\mathcal W}
\newcommand{\cX}{\mathcal X}
\newcommand{\cY}{\mathcal Y}
\newcommand{\px}{p_x}
\newcommand{\py}{p_y}
\newcommand{\sr}{s_r}
\newcommand{\sx}{s_x}
\newcommand{\sy}{s_y}
\newcommand{\diag}{\operatorname{diag}}

\newcommand{\im}{\operatorname{im}}
\DeclareMathOperator{\Bl}{Bl}
\def\point{\mathrm{pt}}

\newtheorem{theorem}{Theorem}[section]
\newtheorem{proposition}[theorem]{Proposition}

\newtheorem{lemma}[theorem]{Lemma}
\theoremstyle{definition}
\newtheorem{remark}[theorem]{Remark}
\theoremstyle{definition}
\newtheorem{definition}[theorem]{Definition}

\title{Rank-Two Attractors and Modularity for a Mirror Bicubic}

\author{Yibing Chen}

\address{Beijing International Center for Mathematical Research, Peking University, 100871, China}

\email{\href{mailto:c_yb@pku.edu.cn}{\nolinkurl{c_yb@pku.edu.cn}}}

\author{Emanuel Scheidegger}

\address{Beijing International Center for Mathematical Research, Peking University, 100871, China}

\email{\href{mailto:esche@bicmr.pku.edu.cn}{\nolinkurl{esche@bicmr.pku.edu.cn}}}

\begin{document}

\begin{abstract}
We study the mirror family $\cX_\lambda$ of the bicubic Calabi--Yau threefold, obtained from a family of complete intersections of two cubics in $\mathbb P^5$.  We show that the four-dimensional weight-three Hodge structure of the fibers over $\lambda^6 = -8$ splits into two rational rank-two pieces of Hodge types $(1,0,0,1)$ and $(0,1,1,0)$; in particular, these points are rank-two attractor points. Using a Hesse pencil fiber-product realization, equivariant point counting, and the Serre--Faltings--Livn\'e method, we identify their semisimplified $\ell$-adic Galois representations with those attached to the weight-four newform $f_{\mathsf{54.4.a.c}}$ and the Tate twist of the weight-two newform $f_{\mathsf{54.2.a.a}}$, where these labels refer to the L-functions and modular forms database.
\end{abstract}

\maketitle

\tableofcontents

\section{Introduction}

The attractor mechanism arose in the study of supersymmetric black holes
in Calabi--Yau compactifications in string theory.  For a smooth projective 
Calabi--Yau threefold, an
attractor point is a distinguished complex structure selected by an
integral class in the middle cohomology.  Moore observed that such points
can exhibit remarkable arithmetic properties and initiated the systematic
study of their relation with complex multiplication and arithmetic
geometry \cite{moore98}.  In the rank-two situation considered here, the
relevant condition is that
$H^{3,0}\oplus H^{0,3}$ be defined over $\mathbb Q$; its intersection
with integral cohomology then gives a lattice of rank two.

A parallel arithmetic question is the modularity of Calabi--Yau
threefolds.  For a rigid Calabi--Yau threefold over $\mathbb Q$,
$H^3$ is two-dimensional and has Hodge type $(1,0,0,1)$, so the
corresponding Galois representation is expected to arise from a
weight-four modular form.  Modularity was established for broad classes
of rigid examples and, after the proof of Serre's conjecture, for rigid
Calabi--Yau threefolds over $\mathbb Q$ in general; see
\cite{dieu03,yui13}.  For a non-rigid threefold the middle cohomology has
larger rank, and modularity can instead appear through lower-dimensional
motivic summands.

Rank-two attractors provide a natural source of such summands. 
If $z_* \in \cM$ is a rank-two attractor point of a family of smooth projective Calabi--Yau threefolds $\cX \to \cM$ over a suitably defined base $\cM$, then the variety $\cX_{z_*}$ is conjectured to be defined over a number field $k \subset \overline{\mQ}$~\cite{Moore:2004fg}. Furthermore, there are three expectations about such points. First, these points are expected to be discrete in $\cM$, see~\cite[Example 5.5]{Klingler:2001}. Second, the varieties $\cX_{z_*}$, if defined over $\mQ$, are expected to be modular in the sense mentioned above. These two expectations were formulated in~\cite{Moore:2004fg}.
Third, attractor points can be viewed as generalization of CM points for elliptic curves. At a CM point, the endomorphism ring of the elliptic curve has additional elements besides the integers. Therefore, it may be expected that the automorphism group of $\cX_{z_*}$ is larger than the automorphism group of a general member of the family. In fact, this need not be the case. Still, it can happen that there exists a family $\cY \to \cM$ such that $k(\cX)$ and $k(\cY)$ have a common finite field extension and it is the automorphism group of $\cY_{z_*}$ that is larger than the automorphism group of a general member $\cY_{z}$. As we will explain below, this is indeed what happens for the family under consideration in this article, and is ultimately the reason for the existence of these arithmetic structures.

The expected presence of these arithmetic structures, if taken for granted, can be used to systematically search for conjectural attractor points. This approach was successfully pursued by Candelas, de la Ossa, Elmi and van Straten in~\cite{candelas2020}, where they found striking persistent quadratic
factorizations of Frobenius polynomials at special fibers of a
one-parameter Calabi--Yau family, leading to rank-two attractor candidates.
Their Hodge-theoretic and arithmetic structure was
subsequently studied by Yang in connection with Deligne's conjecture
\cite{yang2021}.  More recently, Dummigan investigated the modular decomposition
for one of these fibers, obtaining two-dimensional factors associated
with modular forms of weights $4$ and $2$ \cite{dummigan26}.  These
examples show that a rank-two attractor splitting may be reflected
directly in the arithmetic of the middle cohomology.  The result of the
present paper gives another instance of this phenomenon, arising from the
mirror of the bicubic threefold.  We will indeed prove such a weight-four/weight-two
decomposition for one of the previously observed rank-two attractor fibers.

The bicubic Calabi--Yau threefold is a smooth complete intersection of two
cubics in $\mathbb P^5$ given as the special family
\begin{equation}
    \label{eq:Vlambda}
      \begin{aligned}
        \cV_\lambda = \{ x_0^3+x_1^3+x_2^3 - 3\lambda x_3x_4x_5 &= 0,\\
        x_3^3+x_4^3+x_5^3 - 3\lambda x_0x_1x_2 &= 0 \} \subset \mP^5, \quad \lambda^6 \not = 0,1.
      \end{aligned}
\end{equation}
It admits an action by a group $G_{81}$ of order $81$, the largest diagonal symmetry group of $\mP^5$ preserving these equations. Let $\cX_\lambda$ denote the quotient $\cV_\lambda/G_{81}$. Any fiber of this one-parameter family admits a crepant resolution with Hodge numbers $h^{1,1}=73$ and $h^{2,1}=1$ and appears as a mirror to the generic complete intersection of two cubics in $\mP^5$~\cite{libgober1993}. The variation of Hodge structure of this family is governed by the Picard--Fuchs operator
\begin{equation}
\label{eq:mirror-bicubic-PF}
    D_w
=
\theta_w^4
-
9w(3\theta_w+1)^2(3\theta_w+2)^2,
\qquad
\theta_w=w\frac{d}{dw}
\end{equation}
with the parameter
$w=(3\lambda)^{-6}$. This is the hypergeometric Calabi--Yau operator listed as AESZ 4~\cite{AESZ}. Thus AESZ 4 will refer to the Picard--Fuchs equation, whereas the geometric family will be referred to as the mirror bicubic.

Based on the techniques developed in~\cite{candelas2020,Candelas:2021tqt} it was conjectured in~\cite{bksz} that this family admits an attractor point at $w_*=(3\lambda_*)^{-6}=-\frac{1}{5832}$. This is in fact the case:
\begin{theorem}[Main theorem = Theorems~\ref{thm:crepant-resolution-attractors},~\ref{thm:resolved-motivic-modularity}]
\label{thm:attractor-modularity}
Let $\lambda_*=\sqrt{-2}$, and $X=\mathcal V_{\lambda_*}/G_{81}$, 
where $G_{81}$ is an Abelian group of order $81$ introduced in \cref{def:groups-actions}. 
Consider $X^{(0)}$ be the $\mathbb Q$-model of $X$
constructed in Section~\ref{sec:birational-quotients}.
Then the followings are true. 
\begin{enumerate}
\item There exist smooth projective Calabi--Yau threefold
$\widetilde X$ defined over $\mathbb Q$,
together with a smooth projective crepant resolution
$\rho_X:\widetilde X\longrightarrow X^{(0)}$. 
After base change to $\mathbb C$, it is a rank-two attractor in the mirror bicubic family. More precisely, there is a decomposition of rational Hodge structures
\[
H^3\bigl(\widetilde X,\mathbb Q\bigr)=V_{\tilde X,1}\oplus V_{\tilde X,2},
\]
where
\[
V_{\tilde X,1}\otimes_{\mathbb Q}\mathbb C
=H^{3,0}(\tilde X)\oplus H^{0,3}(\tilde X),\qquad
V_{\tilde X,2}\otimes_{\mathbb Q}\mathbb C
=H^{2,1}(\tilde X)\oplus H^{1,2}(\tilde X),
\]
and
\[
\operatorname{rank}_{\mathbb Z}
\bigl(V_{\tilde X,1}\cap H^3(\tilde X,\mathbb Z)\bigr)=2.
\]
In particular, $\lambda_*$ is a rank-two attractor
point of the mirror bicubic family. 

\item Moreover, this splitting is motivic over $\mathbb Q$.
In the category of Andr\'e motives, there are isomorphisms
\[
h^3(\widetilde X)
\simeq \mathsf M_1\oplus\mathsf M_2,
\]
where $\mathsf M_1$ and $\mathsf M_2$ are rank-two pure motives of weight three. For every prime $\ell$, their $\ell$-adic realizations satisfy
\[
H^3_{\mathrm{\acute et}}
\bigl(\widetilde X_{\overline{\mathbb Q}},\mathbb Q_\ell\bigr)
\simeq V_{1,\ell}\oplus V_{2,\ell},
\]
where the two summands are $\operatorname{Gal}(\overline{\mathbb Q}/\mathbb Q)$-representations
associated with $\mathsf M_1$ and $\mathsf M_2$.

\item Writing $f$ and $g$ for the level-$54$ 
newforms of weights $4$ and $2$ given in \cite[\href{https://www.lmfdb.org/ModularForm/GL2/Q/holomorphic/54/4/a/c/}{$f_{\mathsf{54.4.a.c}}$}]{lmfdb}, \cite[\href{https://www.lmfdb.org/ModularForm/GL2/Q/holomorphic/54/2/a/a/}{$f_{\mathsf{54.2.a.a}}$}]{lmfdb}, respectively,
one has
\[
V_{1,\ell}^{\mathrm{ss}}
\simeq \rho_{f,\ell}^{\mathrm{ss}},
\qquad
V_{2,\ell}^{\mathrm{ss}}
\simeq
\bigl(\rho_{g,\ell}(-1)\bigr)^{\mathrm{ss}}
\]
as Galois representations. Consequently, for every prime $p>3$ with $p\ne\ell$,
\[
\det\left(
1-T\operatorname{Frob}_p
\mid
H^3_{\mathrm{\acute et}}
(\widetilde X_{\overline{\mathbb Q}},\mathbb Q_\ell)
\right)=
\bigl(1-a_p(f)T+p^3T^2\bigr)
\bigl(1-p\,a_p(g)T+p^3T^2\bigr)
\]
where $\operatorname{Frob}_p$ denotes the geometric Frobenius. 
\end{enumerate}
\end{theorem}

We now describe the proof. We start with the Hesse pencil of elliptic curves. Its main properties relevant to us and its periods are briefly reviewed in Section~\ref{sec:hesse-pencil}. The variation of Hodge structure underlying the Picard--Fuchs operator $D_w$ in~\eqref{eq:mirror-bicubic-PF} admits a geometric realization in terms of a twisted fiber product $\cW_z$ of two Hesse pencils, with $w=z^3$. 
Taking the Hadamard product of the two second-order Picard--Fuchs operator governing the elliptic periods of either pencil gives the fourth-order operator $D_w$ and identifies the resulting variation with AESZ 4. This construction is the content of Section~\ref{sec:hesse-fiber-products}.

We realize the family $\cW_z$ as a complete intersection in $\mP^2\times\mP^2\times\mP^1$. Quotienting the family $\cW_z$ by the largest diagonal symmetry group of the ambient space preserving $\cW_z$, the group ${G_{27}}\cong\mathbb Z_3^3$, yields a one-parameter family $\cY_z$ of Calabi--Yau threefolds. In Section~\ref{sec:birational-quotients} we show that the families $\cX_\lambda$ and $\cY_z$ are birationally equivalent for $z^3 = (3\lambda)^{-6}$. Moreover, this birational equivalence induces an isomorphism of pure rational Hodge structures of weight three,
\begin{equation}
 H^3(\cX_\lambda,\mathbb Q)\xrightarrow{\sim}H^3(\cY_z,\mathbb Q).
\end{equation}
Furthermore, we show that the singular varieties $\cX_\lambda$ and $\cY_z$ admit crepant smooth projective resolutions $\widetilde \cX_\lambda$ and $\widetilde \cY_z$, respectively, and all varieties are defined over $\mQ$ or have a model over $\mQ$. This is the starting point for the discussion of their arithmetic structures.

The key observation is then that at 
\[
z_*=-\frac{1}{18},
\]
the fiber $W = \cW_{z_*}$ has an extra automorphism subgroup 
$\Gamma\cong\operatorname{SL}(2,\mathbb F_3)$ compared with the generic fiber.  This symmetry allows us to isolate the rank-two attractor summand. In Section~\ref{sec:special-fiber-symmetries} we study the action of
$\Gamma$ on the cohomology of $W$.  Fixed-locus calculations and
the Lefschetz formula determine the character of its action on $H^3(W,\mQ)$.
The invariant part $V_1$ is two-dimensional, and a direct calculation
of the action on the holomorphic three-form shows that its complexification
is exactly $H^{3,0}(W)\oplus H^{0,3}(W)$.  The orthogonal complement
inside $H^3(W)^{{G_{27}}}$ gives $V_2$.  Since $V_1$ is defined
over $\mathbb Q$, its intersection with integral cohomology has rank
two, which proves the attractor statement, for $\cW_{z_*}, \cY_{z_*}$ and $\cX_{\lambda_*}$.

Finally, in Section~\ref{sec:point-counting-modularity} we determine the Frobenius traces on $V_1$ and $V_2$ by equivariant point counting.  Averaging the Lefschetz numbers of $W$ over
${G_{27}}$ gives the trace on the four-dimensional piece corresponding to the mirror bicubic, while averaging over $\Gamma$ gives the trace on $V_1$.
For the non-diagonal $\Gamma$-action, the required semilinear fixed
points are obtained from the Lang--Steinberg theorem. The resulting traces agree with the Fourier coefficients of
$f_{\mathsf{54.4.a.c}}$ and with
$p\,a_p(f_{\mathsf{54.2.a.a}})$ given in~\cite{lmfdb}. The
Serre--Faltings--Livn\'e criterion, in the finite form used by Sch\"utt,
then proves the modularity statement.

It turns out that the same construction and steps can be applied to other cases. For example, there exists another fiber product of two twisted Hesse pencils. The Hadamard product construction then leads to the Picard--Fuchs operator denoted AESZ 165 in~\cite{AESZ}. For this operator there is a conjectural rank-two attractor point listed in~\cite{TRCP22}. We plan to report on this case in~\cite{ChenScheidegger2}. 


\textbf{Acknowledgments:}
We have greatly benefitted from discussions with Matt Kerr, Aimeric Malter, Joseph McGovern, Yuze Sun and Liang Xiao. Parts of the birational equivalence in Theorem~\ref{thm:quotient-H3-comparison} have independently been proven by Aimeric Malter and we would like to thank him for informing us. The work of the authors is partially supported by NSFC grant 12471059. 

\textbf{AI statement:}
The main mathematical ideas and the initial draft of this paper were developed by the human authors before the rise of AI in math. AI tools were used to assist with writing PARI/GP code, checking computational results, and reviewing technical details for consistency and accuracy. The authors independently reviewed the AI--assisted work and take full responsibility for the content and correctness of the paper. 

\section{Hesse pencil and its periods}
\label{sec:hesse-pencil}

In this section, we review the basic geometry of the Hesse pencil of elliptic curves. We will discuss the singular fibers, torsion sections, and automorphism groups. Then we derive the Picard--Fuchs equation of the periods of the holomorphic one-form and express the holomorphic solution in terms of a hypergeometric function. These results provide backgrounds for future constructions of the Calabi--Yau threefolds. 

\begin{definition}[The Hesse pencil]

Let $\xi_3=\exp(2\pi i/3)=(-1+\sqrt{-3})/2$ be a third root of unity. 
The Hesse pencil is the elliptic surface $\mathcal E$ together with the projection map $\operatorname{pr}$. 
\begin{align}
  \mathcal E&:=\{((x_0:x_1:x_2),(t_0:t_1))\in\mathbb P^2\times\mathbb P^1:
 t_0(x_0^3+x_1^3+x_2^3)=t_1x_0x_1x_2\}, 
  \label{eq:hesse-pencil-family}\\
  \nonumber
\operatorname{pr}&:\mathcal E\to\mathbb P^1,\qquad
\operatorname{pr}((x_0:x_1:x_2),(t_0:t_1)):=(t_0:t_1).
\end{align}

For simplicity we set $t_1=1$, and take as parameter $t=t_0/t_1$, so that the fiber is the plane cubic
\[
\mathcal E_t:=\{(x_0:x_1:x_2)\in\mathbb P^2: t(x_0^3+x_1^3+x_2^3)=x_0x_1x_2\}. 
\]
We also fix the base point $O=(0:1:-1)$ on each fiber. 
\end{definition}

Now we describe the geometry in the following proposition, these are the standard results of the Hesse pencil, see \cite[Section~2]{artebani-dolgachev09}. The formulas below are written in our choice of parameter and coordinates. 

\begin{proposition}

For $t_1^3\ne 27t^3, t_0^3\ne 0$, the $j$-invariant of the fiber elliptic curve $\mathcal E_{t_0/t_1}$ is

\begin{align}
 j(\mathcal E_{t_0/t_1})
 =\frac{t_1^3(t_1^3+216t_0^3)^3}{t_0^3(t_1^3-27t_0^3)^3}.        \label{eq:hesse-j-invariant}
\end{align}

The singular fibers occur at
\[
(t_0:t_1)=(0:1),\quad (1:3),\quad (\xi_3^2:3),\quad (\xi_3:3),
\]

and are given by
\begin{align*}
\mathcal E_{0}:&\quad x_0x_1x_2=0,\\
\mathcal E_{1/3}:&\quad
(x_0+x_1+x_2)(x_0+\xi_3x_1+\xi_3^2x_2)(x_0+\xi_3^2x_1+\xi_3x_2)=0,\\
\mathcal E_{\xi_3^2/3}:&\quad
(\xi_3x_0+x_1+x_2)(\xi_3x_0+\xi_3x_1+\xi_3^2x_2)(\xi_3x_0+\xi_3^2x_1+\xi_3x_2)=0,\\
\mathcal E_{\xi_3/3}:&\quad
(\xi_3^2x_0+x_1+x_2)(\xi_3^2x_0+\xi_3x_1+\xi_3^2x_2)(\xi_3^2x_0+\xi_3^2x_1+\xi_3x_2)=0.
\end{align*}

Each singular fiber is of Kodaira type $I_3$.  Moreover $\mathcal E$ is smooth
and the nine sections
\begin{align}
 (\xi_3^k:-1:0),\qquad (0:\xi_3^k:-1),\qquad (-1:0:\xi_3^k),
 \qquad k=0,1,2,                                          \label{eq:hesse-torsion-sections}
\end{align}
realize the full $3$-torsion subgroup of each smooth fiber.  In particular,
$\mathcal E$ is the universal elliptic curve over $X(3)$.
\end{proposition}

We now record the automorphisms of the Hesse pencil that will be used later.
They will reappear in Section~\ref{sec:middle-cohomology-actions} in the study of automorphism groups acting on the cohomology of the associated Calabi--Yau threefolds.

\begin{definition}[Hessian group]
The Hessian group is the group of automorphisms of the Hesse pencil 
$\mathcal E$ generated by translations by the nine sections in
\eqref{eq:hesse-torsion-sections} and by the automorphisms induced by linear transformations
of $\mathbb P^2$ preserving the Hesse pencil.  It is isomorphic to
\[
  \mathrm{Hes}:=\operatorname{ASL}(2,\mathbb F_3)
  \cong \mathbb Z_3^2\rtimes \operatorname{SL}(2,\mathbb F_3),
\]
and has order $9\cdot 24=216$.
\end{definition}

\begin{proposition}[Automorphisms of the Hesse pencil]
\label{prop:hesse-automorphisms}
The following assertions hold.
\begin{enumerate}
\item Fiberwise, the translations by the $3$-torsion sections together with the
elliptic involution generate an order-$18$ subgroup
\[
  \mathcal E_t[3]\rtimes\{\pm1\}
  \cong \mathbb Z_3^2\rtimes \mathbb Z_2.
\]
It may be realized inside $\operatorname{PGL}(3,\mathbb C)$ by the permutations
of $x_0,x_1,x_2$ and by the diagonal transformation
$(x_0:x_1:x_2)\mapsto (x_0:\xi_3x_1:\xi_3^2x_2)$.

\item The two transformations
\begin{align}
\label{eq:hesse-phi-action}
&\phi:\quad
\begin{cases}(t_0:t_1)\mapsto(\xi_3t_0:t_1),\\
(x_0:x_1:x_2)\mapsto(\xi_3x_0:x_1:x_2),\end{cases} \\
\label{eq:hesse-psi-action}
&\psi:\quad
\begin{cases}(t_0:t_1)\mapsto(t_1-3t_0:3t_1+18t_0),\\
(x_0:x_1:x_2)\mapsto(x_0+x_1+x_2:x_0+\xi_3x_1+\xi_3^2x_2:x_0+\xi_3^2x_1+\xi_3x_2)\end{cases}
\end{align}
are automorphisms of the total space $\mathcal E$ that 
maps the base point section $O\times\mathbb P^1$ to itself. 
Their induced action on the base generates a group 
isomorphic to $A_4\cong \operatorname{PSL}(2,\mathbb F_3)$.
This is the Galois group of the degree-$12$ map
\[
  j:X(3)\longrightarrow X(1).
\]

\item The full group generated by the fiberwise automorphisms and the two
transformations above is the Hessian group
$\mathrm{Hes}=\mathbb Z_3^2\rtimes\operatorname{SL}(2,\mathbb F_3)$. 
And the subgroup $\Gamma\subset\mathrm{Hes}$ generated by $\phi,\psi$ is isomorphic to $\operatorname{SL}(2,\mathbb F_3)$. 
\end{enumerate}

\end{proposition}

\begin{proof}

The first statement follows from the fact that the flex points of a smooth
Hesse cubic are exactly the $3$-torsion points.  The permutations of the three
coordinates and the diagonal map preserve the Hesse pencil and realize the
translations by $3$-torsion together with inversion.

For the second statement one substitutes \eqref{eq:hesse-phi-action} and \eqref{eq:hesse-psi-action} into
$t_0(x_0^3+x_1^3+x_2^3)=t_1x_0x_1x_2$.  In both cases the transformed equation is again the
same equation with the transformed base point. Set $t_1=1$, the induced transformations on
$t=t_0/t_1$ are fractional linear transformations
\[
  t\mapsto \xi_3t,
  \qquad
  t\mapsto \frac{1-3t}{3+18t}.
\]

These maps preserve the four cusps
$\{0,\frac{1}{3},\frac{\xi_3}{3},\frac{\xi_3^2}{3}\}$ and act alternatingly. 
Hence it is easy to check that they generate the alternating group
$A_4$. Since the map $X(3)\to X(1)$ has degree
$|\operatorname{PSL}(2,\mathbb F_3)|=|A_4|=12$, this group is its Galois group.

The last assertion is the standard description of the Hessian group as the
semidirect product of the translation group by the natural
$\operatorname{SL}(2,\mathbb Z)/\Gamma(3)\cong\operatorname{SL}(2,\mathbb F_3)$-action on $\mathcal E_t[3]$, see~\cite[Section 4]{artebani-dolgachev09}. Now we turn to the automorphism subgroup $\Gamma:=\langle\phi,\psi\rangle$; its generators satisfy the structural relations:
\begin{align}
\label{eq:gamma-relations}
\phi^3 = \psi^4 = (\phi\psi)^3 = [\phi, \psi^2] = 1.
\end{align}

Standard computations in group theory demonstrate that the free group modulo these relations is isomorphic to the order $24$ group $\operatorname{SL}(2,\mathbb F_3)$. It is easy to check that $\Gamma$ is a quotient of $\operatorname{SL}(2,\mathbb F_3)$ with $24$ distinct elements, hence $\Gamma\cong\operatorname{SL}(2,\mathbb F_3)$. 
\end{proof}

On each smooth fiber $\mathcal E_t$, we define the
holomorphic one-form $\omega_t$ as the Poincaré residue
along $\mathcal E_t$ of a given rational two-form on $\mathbb P^2$:
\begin{align}
\label{eq:hesse-holomorphic-form}
\omega_t:=&\operatorname{Res}
\frac{x_0\,dx_1\,dx_2+x_1\,dx_2\,dx_0+x_2\,dx_0\,dx_1}{x_0x_1x_2-t(x_0^3+x_1^3+x_2^3)}\\
=&\operatorname{Res}\left(\frac{dx_1}{x_1}\frac{dx_2}{x_2}
 +\frac{dx_2}{x_2}\frac{dx_0}{x_0}+\frac{dx_0}{x_0}\frac{dx_1}{x_1}\right)
\left(1-t\cdot\frac{x_0^3+x_1^3+x_2^3}{x_0x_1x_2}\right)^{-1}. 
\nonumber
\end{align}

\begin{proposition}[Picard--Fuchs equation]
Let $\theta_t=t\frac{d}{dt}$. On the smooth locus $t\ne0$ and $27t^3\ne 1$, the Picard--Fuchs operator satisfied by the holomorphic one-form $\omega_t$ is
\begin{align}
L_t:=\theta_t^2-27t^3(\theta_t+1)(\theta_t+2). \label{eq:hesse-PF-operator}
\end{align}
More precisely, for every locally constant family of cycles $\gamma_t\in H_1(\mathcal E_t,\mathbb Z)$, the period $\varpi_\gamma(t):=\int_{\gamma_t}\omega_t$ satisfies
\[
L_t\varpi_\gamma(t)=0. 
\]
\end{proposition}

\begin{proof}
With $\psi=(3t)^{-1}$, our form $\omega_t$ equals $-3$ times the one-form in \cite[(1.1)]{zhou17}, whose periods are annihilated by $(\theta_\psi-2)(\theta_\psi-1)-\psi^3\theta_\psi^2$ by \cite[(1.2)]{zhou17}. 
Since $\theta_\psi=-\theta_t$, multiplying this operator by $-\psi^{-3}=-27t^3$ gives $\theta_t^2-27t^3(\theta_t+1)(\theta_t+2)=L_t$, proving the assertion.
\end{proof}

In a neighbourhood of $t=0$, the unique holomorphic solution of $L_t$ up to a scalar is
\begin{align*}
\varpi(t)&={}_2F_1\!\left({{\frac{1}{3},\frac{2}{3}}\atop {1}};(3t)^3\right)
 =\sum_{n=0}^\infty\frac{(3n)!}{(n!)^3}t^{3n}   \\
&=1+6t^3+90t^6+1680t^9+34650t^{12}+756756t^{15}+17153136t^{18}+O(t^{21}), 
\end{align*}
this is actually the generating function of the constant terms of
$\left(\frac{x_0^3+x_1^3+x_2^3}{x_0x_1x_2}\right)^n$. 

\section{Fiber products of twisted Hesse pencils}
\label{sec:hesse-fiber-products}

In this section, we construct a family of Calabi--Yau threefolds as Schoen fiber products of two twisted Hesse pencils and study its basic properties. We then construct the holomorphic three-form and compute its periods using the Hadamard product of the elliptic periods. This leads to a fourth-order Picard--Fuchs operator, which we identify with the AESZ 4 operator associated with the mirror bicubic. 

\begin{definition}
For any rational elliptic surfaces $\pi_1: \mathcal E^{(1)}\to\mathbb P^1$ and $\pi_2:\mathcal E^{(2)}\to\mathbb P^1$ with sections $O_1,O_2$, respectively, we define the Schoen fiber product $X$ of $\mathcal E^{(1)}$ and $\mathcal E^{(2)}$ as
\begin{align}
X:=\mathcal E^{(1)}\times_{\mathbb P^1} \mathcal E^{(2)}.
\label{eq:schoen-fiber-product}
\end{align}
\end{definition}

These fiber product threefolds are in general not smooth. Singular points occur precisely over those $t\in\mathbb P^1$ for which both $\mathcal E^{(1)}_t$ and
$\mathcal E^{(2)}_{t}$ are singular. If the two singular fibers are of
Kodaira types $I_a$ and $I_b$, then their product contains $ab$ ordinary
double points. In our cases only semistable fibers of type $I_n$ occur, then after small resolution, this is exactly Schoen's construction of Calabi--Yau threefolds. 

\begin{theorem}[Schoen]
\label{thm:schoen-small-resolutions}
After performing small projective resolutions of the ordinary double points of $X$, the resolved variety $\widetilde X$ is a smooth projective Calabi--Yau
threefold. 
Moreover, $H^2(\widetilde X,\mathbb Q)$ is generated by the
classes of the following divisors, and
$H^4(\widetilde X,\mathbb Q)$ is generated by their
intersections with an ample divisor:

\begin{enumerate}
\item a fiber of the fibration
$\mathcal E^{(1)}\times_{\mathbb P^1} \mathcal E^{(2)}\to\mathbb P^1$;
\item all but one irreducible components of each reducible fiber;
\item the two product divisors
$\mathcal E^{(1)}\times_{\mathbb P^1}O_2$ and $O_1\times_{\mathbb P^1} \mathcal E^{(2)}$; 
\item Mordell Weil divisors: if $\{P_1,\cdots,P_r\}$ is a basis of the free part of $\operatorname{MW}(\mathcal E^{(1)})$, $\{Q_1,\cdots,Q_s\}$ is a basis of the free part of $\operatorname{MW}(\mathcal E^{(2)})$, then the divisors are
$\mathcal E^{(1)}\times_{\mathbb P^1}Q_j$ and $P_i\times_{\mathbb P^1} \mathcal E^{(2)}$; 
\item closure of graphs of isogenies between the two generic elliptic fibers, if such isogenies exist.
\end{enumerate}
\end{theorem}

\begin{proof}
Since $X$ is the inverse image of the diagonal in $\mathbb P^1\times\mathbb P^1$ under $\pi_1\times\pi_2$, it is a Cartier divisor in the smooth fourfold $\mathcal E^{(1)}\times\mathcal E^{(2)}$.
Thus $X$ is Gorenstein, and adjunction gives
\[
\omega_X\cong f^*\left(\mathcal O_{\mathbb P^1}(-1)\otimes\mathcal O_{\mathbb P^1}(-1)\otimes\omega_{\mathbb P^1}^{-1}\right)\cong\mathcal O_X,
\]
where $f:X\to\mathbb P^1$ is the natural projection.
The ordinary double points are rational singularities, and any small resolution $\rho:\widetilde X\to X$ is crepant, so $\omega_{\widetilde X}\cong\mathcal O_{\widetilde X}$ and $R\rho_*\mathcal O_{\widetilde X}
\simeq\mathcal O_X$.
Moreover, by flat base change and the relative K\"unneth formula,
\[
R^0f_*\mathcal O_X\cong\mathcal O_{\mathbb P^1},
\quad
R^1f_*\mathcal O_X\cong
\mathcal O_{\mathbb P^1}(-1)^{\oplus 2},
\quad
R^2f_*\mathcal O_X\cong
\mathcal O_{\mathbb P^1}(-2).
\]
The Leray spectral sequence for $f$ therefore gives
$H^1(X,\mathcal O_X)=H^2(X,\mathcal O_X)=0$,
and hence the same vanishing holds for
$\widetilde X$. Thus $\widetilde X$ is a smooth projective Calabi--Yau threefold. For more details see \cite{sch88}. 

Let $K=k(\mathbb P^1)\cong\mathbb C(t)$ and let $E_1,E_2$ be the generic elliptic fibers over $K$ with sections. The Picard group of $E_1\times_{\operatorname{Spec}K} E_2$ is generated by the two factors, the rational $K$-points of $E_1$ and $E_2$, and the graphs of homomorphisms between them, up to torsion. These give the divisors in (3)--(5). 
The kernel of restriction to the generic fiber is generated by vertical divisors, namely
a fiber and all but one component of each reducible fiber. These give (1)--(2). 
Consequently, the listed divisors generate $\operatorname{NS}(\widetilde X)\otimes\mathbb Q$, as in Schoen's description \cite{sch88}. 

Since $H^{2,0}(\widetilde X)=0$, the Lefschetz
$(1,1)$-theorem identifies this space with
$H^2(\widetilde X,\mathbb Q)$.
Finally, for any ample divisor $D$, the hard Lefschetz
theorem gives an isomorphism
\[
D\cup-:
H^2(\widetilde X,\mathbb Q)
\xrightarrow{\sim}
H^4(\widetilde X,\mathbb Q).
\]
Hence $H^4$ is generated by the classes
obtained by intersecting $D$ with the listed divisors.
\end{proof}

\begin{definition}
\label{def:twisted-hesse-product}
For $z\notin \Delta:=
\left\{0,\frac1{3^2},\frac{\xi_3}{3^2},
\frac{\xi_3^2}{3^2},\infty\right\}$, let $\iota_z: \mP^1 \to \mP^1$ denote the map $t\mapsto \frac{z}{t}$. We define the twisted Schoen fiber product of the Hesse pencil $\cE$ in~\eqref{eq:hesse-pencil-family}:
\[
\mathcal W_z:=\mathcal E\times_{\mathbb P^1} \iota_z^*\mathcal E. 
\]
\end{definition}
In terms of the realization of the rational elliptic surface $\cE$ as a hypersurface in $\mP^2 \times \mP^1$ we have
\begin{equation}
  \label{eq:Wz_CICY}
  \cW_z =\{t_0(x_0^3+x_1^3+x_2^3)=t_1x_0x_1x_2,\quad zt_1(y_0^3+y_1^3+y_2^3)=t_0y_0y_1y_2\} \subset \mP^2 \times \mP^2 \times \mP^1
\end{equation}

\begin{proposition}
For $z\notin \Delta$, there are no common singular fibers of $\mathcal E_t$ and $\mathcal E_{z/t}$. Hence $\mathcal W_z$ is a smooth projective Calabi--Yau threefold and
\[
h^{1,1}\left(\mathcal W_z\right)=h^{2,1}\left(\mathcal W_z\right)=19,\qquad \chi\left(\mathcal W_z\right)=0.
\]
\end{proposition}

\begin{proof}
The singular fibers $\mathcal E_t$ are at $t \in \{0,\frac13,\frac{\xi_3}{3},\frac{\xi_3^2}{3}\}$. For $z\neq0,\infty$, the singular fibers $\mathcal E_{z/t}$, regarded as an elliptic surface over the $t$-line, are at $t\in\{\infty,3z,3\xi_3z,3\xi_3^2z\}$. Thus a common singular fiber can occur only when
\[
t\in\left\{\frac13,\frac{\xi_3}{3},\frac{\xi_3^2}{3}\right\},\qquad \frac{z}{t}\in\left\{\frac13,\frac{\xi_3}{3},\frac{\xi_3^2}{3}\right\},
\]
or equivalently when $z\in\{\frac1{3^2},\frac{\xi_3}{3^2},\frac{\xi_3^2}{3^2}\}$. Together with $z=0,\infty$, these are precisely the values in $\Delta$. Since the fiber product of two smooth surfaces over a smooth curve is singular at a point only if both projections are critical there, it follows that $\mathcal E\times_{\mathbb P^1} \iota_z^*\mathcal E$ is smooth for $z\notin\Delta$ and hence $\mathcal W_z$ is a smooth projective Calabi--Yau threefold. 

We next compute $h^{1,1}$ by \cref{thm:schoen-small-resolutions}. Each of the two rational elliptic surfaces has four fibers of type $I_3$, and the two sets of singular values are disjoint. 
They are shown as follows. 
\begin{align*}
\begin{array}{ccccccccc}
t & 0 & 1/3 & \xi_3/3 & \xi_3^2/3 & 3z & 3\xi_3 z & 3\xi_3^2 z & \infty \\
\hline
\mathcal E_t & I_3 & I_3 & I_3 & I_3 & - & - & - & -\\
\mathcal E_{z/t} & - & - & - & - & I_3 & I_3 & I_3 & I_3
\end{array}
\end{align*}

By the Shioda--Tate formula and $\mathcal E\cong\operatorname{Bl}_9\mathbb P^2$, 
\[
9+1=\rho(\mathcal E)=2+4(3-1)+\operatorname{rk}\operatorname{MW}(\mathcal E/\mathbb P^1),
\]
so the Mordell--Weil rank of either elliptic surface is zero. Moreover, if $E_1$ and $E_2$ denote their generic elliptic fibers over $\mathbb C(t)$, then $\operatorname{Hom}_{\mathbb C(t)}(E_1,E_2)=0$. A nonzero homomorphism is an isogeny. Note that isogenous elliptic curves have the same places of good reduction and the bad reduction loci of $E_1$ and $E_2$ are disjoint and nonempty. Consequently
\[
h^{1,1}(\mathcal W_z)=1+8(3-1)+2+0+0=19.
\]

Finally, consider the fibration $\mathcal W_z\to\mathbb P^1$. Over every point of the base at least one of the two elliptic fibers is smooth, and a smooth elliptic curve has Euler characteristic zero. Hence every fiber of $\mathcal W_z\to\mathbb P^1$ has Euler characteristic zero, by additivity of the topological Euler characteristic,
$\chi(\mathcal W_z)=0$. For a Calabi--Yau threefold one has $\chi(\mathcal W_z)=2\bigl(h^{1,1}(\mathcal W_z)-h^{2,1}(\mathcal W_z)\bigr)$. Therefore $h^{2,1}(\mathcal W_z)=19$, as claimed.
\end{proof}

\begin{proposition}
\label{prop:holomorphic-period-Wz}
For $z\notin\Delta$, the form
\begin{align}
\label{eq:holomorphic-three-form-Wz}
\Omega_z:=\omega_t\wedge\omega_{z/t}\wedge\frac{dt}{t}
\end{align}
defined over the locus where both elliptic fibers are smooth extends uniquely to a nowhere vanishing holomorphic three-form on $\mathcal W_z$. In particular, it spans $H^{3,0}(\mathcal W_z)$. Moreover, its periods over suitable product cycles are Hadamard products of the periods of the two elliptic curves. If
\[
\varpi(t)={}_2F_1\!\left({\frac13,\frac23\atop 1};(3t)^3\right)
 =\sum_{n=0}^\infty\frac{(3n)!}{n!^3}t^{3n} 
\]
denotes the normalized holomorphic period of the Hesse pencil near $t=0$, then the corresponding holomorphic period of $\Omega_z$ near $z=0$ is
\[
(\varpi\star\varpi)(z)
=\sum_{n=0}^\infty\frac{(3n)!^2}{(n!)^6}z^{3n}
={}_4F_3\left(\begin{matrix}\frac13,\frac13,\frac23,\frac23\\1,1,1\end{matrix};729z^3\right).
\]
\end{proposition}

\begin{proof}
Away from the singular fibers, $\Omega_z=\omega_t\wedge\omega_{z/t}\wedge dt/t$ is clearly a nonzero holomorphic three-form. We check that it extends across the singular fibers. Let $a\neq0,\infty$ be a singular value of one of the two elliptic factors. Since every singular fiber is of type $I_3$, near each node there are local coordinates $u,v$ on the elliptic surface and a local coordinate $q$ on the base such that $uv=q$, while the relative holomorphic one-form is, up to a holomorphic unit $du/u$. If $\mathcal E_t$ is singular at $t=a$, then $q=t-a$ and the other elliptic factor is smooth since $z\notin\Delta$; hence
\[
\Omega_z\sim\frac{du}{u}\wedge\omega_{z/t}\wedge dq=\;du\wedge\omega_{z/t}\wedge dv,
\]
up to a nonzero holomorphic factor. Thus $\Omega_z$ extends holomorphically and does not vanish. The same argument applies when $\mathcal E_{z/t}$ is singular, since $z/t$ is a local coordinate on the base up to a nonzero holomorphic factor.

It remains to consider $t=0,\infty$. For a parameter $u$ near infinity, putting $q=1/u$ in the residue expression~\eqref{eq:hesse-holomorphic-form} gives
\[
\omega_u=-q\,\widehat\omega_q,\qquad
\widehat\omega_q:=\operatorname{Res}\frac{x\,dy\,dz+y\,dz\,dx+z\,dx\,dy}{x^3+y^3+z^3-qxyz},
\]
where $\widehat\omega_q$ is holomorphic and nonzero near $q=0$. Hence, as $t\to0$,
\[
\omega_{z/t}=-\frac{t}{z}\widehat\omega_{t/z},\qquad
\Omega_z=-\frac1z\,\omega_t\wedge\widehat\omega_{t/z}\wedge dt.
\]
Near a node of the $I_3$-fiber $\mathcal E_0$, write $uv=t$. Since $\omega_t\sim du/u$, we obtain
\[
\omega_t\wedge dt\sim\frac{du}{u}\wedge d(uv)=du\wedge dv,
\]
so $\Omega_z$ extends holomorphically and nonvanishingly across $t=0$. Finally, at $t=\infty$ put $q=1/t$. Then $\omega_t=-q\widehat\omega_q$ and $dt/t=-dq/q$, and therefore
\[
\Omega_z=\widehat\omega_q\wedge\omega_{zq}\wedge dq.
\]
Here $\mathcal E_t$ is smooth at $q=0$, whereas $\mathcal E_{zq}$ has an $I_3$-fiber; applying again the local model $uv=q$ shows that $\omega_{zq}\wedge dq$ is holomorphic and nonzero. Thus $\Omega_z$ extends to a global nowhere vanishing holomorphic three-form on $\mathcal W_z$.

We finally compute its periods. For $|z|$ sufficiently small choose $r$ such that $3|z|<r<1/3$, and let $\gamma_t\subset\mathcal E_t$ and $\gamma_{z/t}\subset\mathcal E_{z/t}$ be locally constant families of one-cycles along the circle $|t|=r$. Fibering $\gamma_t\times\gamma_{z/t}$ over this circle gives a three-cycle $\Gamma_z\subset\mathcal W_z$. If the cycles are normalized so that $\int_{\gamma_t}\omega_t=\varpi(t)$, then Fubini's theorem gives
\[
\frac{1}{2\pi i}\int_{\Gamma_z}\Omega_z
=\frac{1}{2\pi i}\oint_{|t|=r}\varpi(t)\varpi(z/t)\frac{dt}{t}.
\]
Writing $\varpi(t)=\sum_{n\geq0}c_n t^{3n}$, where $c_n=(3n)!/(n!)^3$, one obtains
\[
\varpi(t)\varpi(z/t)
=\sum_{m,n\geq0}c_mc_nz^{3n}t^{3m-3n}.
\]
The contour integral extracts the constant term in $t$, hence only the terms with $m=n$ survive, and therefore
\[
\frac{1}{2\pi i}\int_{\Gamma_z}\Omega_z
=\sum_{n\geq0}c_n^2z^{3n}
=\sum_{n\geq0}\frac{(3n)!^2}{(n!)^6}z^{3n}
=(\varpi\star\varpi)(z)
\]
is the holomorphic period of $\mathcal W_z$. 
\end{proof}

By this Hadamard product construction, we can compute the Picard--Fuchs
equation for three-form of $\mathcal W_z$ as the Hadamard product
$D_z=L_z\star L_z$.  
Writing $\theta_z=z\frac{d}{dz}$, the coefficients
$a_n=(3n)!^2/(n!)^6$ satisfy the recursive relation
\[
(3n)^4a_n
=729(3n-2)^2(3n-1)^2a_{n-1}.
\]
It follows that $(\varpi\star\varpi)(z)$ is annihilated by the fourth-order operator
\[
D_z
=\theta_z^4
-729z^3(\theta_z+1)^2(\theta_z+2)^2,
\qquad
D_z {}_4F_3\left(\begin{matrix}\frac13,\frac13,\frac23,\frac23\\1,1,1\end{matrix};729z^3\right)=0.
\]
Equivalently, in the parameter $w=z^3$, so that
$\theta_z=3\theta_w$, this becomes
\[
D_w
=\theta_w^4
-9w(3\theta_w+1)^2(3\theta_w+2)^2,
\qquad
\theta_w=w\frac{d}{dw}.
\]
This is the Picard--Fuchs operator~\eqref{eq:mirror-bicubic-PF} given in the Introduction.

\section{The quotients and birational equivalences}
\label{sec:birational-quotients}

We now have two families of projective Calabi--Yau threefolds, $\cX_\lambda$ and $\cW_z$, which each have a period integral that is annihilated by the Picard--Fuchs operator~\eqref{eq:mirror-bicubic-PF}, up to a change of variables. 
In this section we show that the quotient families $\cX_\lambda=\cV_\lambda/G_{81}$ with $\cV_\lambda$ given in~\eqref{eq:Vlambda} and $\cY_z=\cW_z/G_{27}$ with $\cW_z$ in~\eqref{eq:Wz_CICY} are birational if their parameters are related by $z^3 = (3\lambda)^{-6}$. Then we show that this birational equivalence induces isomorphisms of rational Hodge structures on their middle cohomologies, in particular these are pure despite the fact that the varieties are singular.

\subsection{The birational equivalence}
\label{sec:Birat_equiv}

We start with fixing some notation. We set $x=(x_0,x_1,x_2)$ and $y =(y_0,y_1,y_2)$ for the homogeneous coordinates on the two factors of $\mP^2$ and $(t_0,t_1)$ for $\mP^1$ of the ambient space $\mP^2\times\mP^2\times\mP^1$ of the family $\cW_z$ in~\eqref{eq:Wz_CICY}. As before, we will write $t=t_0/t_1$.
For the homogeneous coordinates of the ambient space $\mP^5$ of the family $\cV_\lambda$ we write $(x,y) = (x_0,x_1,x_2,y_0,y_1,y_2)$. Put
\[
    \px:=x_0x_1x_2, \quad \py:=y_0y_1y_2, \quad
    \sx:=\sum_{i=0}^2x_i^3, \quad \sy:=\sum_{j=0}^2y_j^3.
\]
The variables $\px, \py, \sx, \sy$ will later serve as local coordinates. For completeness, we repeat the equations of the two families here
\begin{align}
  \cV_\lambda&:=\{\sx-3\lambda\py=\sy-3\lambda\px=0\}\subset\mP^5,
 \label{eq:Vlambda-equations}\\
  \cW_z&:=\{t_0\sx - t_1\px=zt_1\sy - t_0\py=0\}
 \subset\mP^2\times\mP^2\times\mP^1, 
 \label{eq:Wz-equations}    
\end{align}
for the choice of parameters
\begin{equation} 
 \label{eq:generic-parameters}
 \lambda\in\mC^*,\quad \lambda^6\ne1;
 \qquad z\in\mC^*,\quad (9z)^3\ne1.
\end{equation}

\begin{definition}
\label{def:groups-actions}
Define the group
\[
{G_{81}}:=
\frac{\left\{(a_0,a_1,a_2,b_0,b_1,b_2)\in\mathbb Z_9^6:
\begin{matrix} 3a_0=3a_1=3a_2=b_0+b_1+b_2,\\ 3b_0=3b_1=3b_2=a_0+a_1+a_2 \end{matrix} \right\}}
{\langle(1,1,1,1,1,1)\rangle}
\cong\mathbb Z_3^2\times\mathbb Z_9,
\]
which acts on $\mP^5$ by 
\[
\begin{aligned}
&(a_0,a_1,a_2,b_0,b_1,b_2):\\
&\quad
(x_0,x_1,x_2,y_0,y_1,y_2)
\longmapsto
\bigl(
\xi_9^{a_0}x_0,\,
\xi_9^{a_1}x_1,\,
\xi_9^{a_2}x_2,\,
\xi_9^{b_0}y_0,\,
\xi_9^{b_1}y_1,\,
\xi_9^{b_2}y_2
\bigr), 
\end{aligned}
\]
Here $\xi_9=\exp(2\pi i/9)$ is a ninth root of unity. It is easy to see that the action is faithful and projective on $\cV_\lambda$, see also~\cite[Section 3.2]{malter2024}. We can rewrite it as
\begin{equation}
\label{eq:G81-action}
\begin{split}
 g(\alpha,\beta,\gamma)=\diag\bigl(\xi_9^{3\alpha+\gamma},
 \xi_9^{-3\alpha+4\gamma},\xi_9^\gamma, \xi_9^{-3\beta-\gamma},\xi_9^{3\beta-4\gamma},\xi_9^{-\gamma}\bigr),
\end{split}
\end{equation}
where $\alpha,\beta\in\mZ_3$ and $\gamma\in\mZ_9$. 

Similarly, we define the group
\[
{G_{27}}:=
\frac{\{(a_0,a_1,a_2,b_0,b_1,b_2)\in\mathbb Z_3^6:
a_0+a_1+a_2+b_0+b_1+b_2=0\}}
{\langle(1,1,1,0,0,0),(0,0,0,1,1,1)\rangle}
\cong\mathbb Z_3^3,
\]
which acts on $\mP^2\times\mP^2\times\mP^1$ by
\[
\begin{aligned}
&(a_0,a_1,a_2,b_0,b_1,b_2):\\
&\quad
(x_0,x_1,x_2,y_0,y_1,y_2,t)
\longmapsto
\bigl(
\xi_3^{a_0}x_0,\,
\xi_3^{a_1}x_1,\,
\xi_3^{a_2}x_2,\,
\xi_3^{b_0}y_0,\,
\xi_3^{b_1}y_1,\,
\xi_3^{b_2}y_2,\,
\xi_3^{a_0+a_1+a_2}t
\bigr).
\end{aligned}
\]
\end{definition}
This group action is also faithful and projective on $\cW_z$ and can be rewritten as
\begin{equation}
\label{eq:G27-action}
\begin{split}
 h(\alpha,\beta,\gamma):(x;y;t_0:t_1)\longmapsto
 \bigl(\xi_3^\gamma x_0:\xi_3^\alpha x_1:\xi_3^{-\alpha}x_2;\xi_3^{-\gamma}y_0:\xi_3^\beta y_1:\xi_3^{-\beta}y_2;
 \xi_3^\gamma t_0:t_1\bigr),
\end{split}
\end{equation}
with $\alpha,\beta,\gamma\in\mZ_3$. 
\begin{definition}
We define the quotient varieties
\[
\begin{aligned}
  \cX_\lambda &:= \cV_\lambda/G_{81}, &  \cY_z &:= \cW_z/G_{27}   
\end{aligned}
\]
\end{definition}

\begin{remark}
\begin{enumerate}
   \label{rem:toric}
   \item The actions of $G_{81}$ on $\cV_\lambda$ and $G_{27}$ on $\cW_z$ both have fixed points, hence both $\cX_\lambda$ and $\cY_z$ are singular. The varieties $\cX_\lambda$ have six rational curves of $A_2$--singularities and six singular points of the form $\mA^3/(\mZ_9 \times \mZ_3)$, with three singular curves meeting at each such point, see~\cite{libgober1993}. Similarly, it is straightforward to see that the varieties $\cY_z$ have 24 rational curves of $A_2$--singularities and 18 singular points of the form $\mA^3/(\mZ_3)^2$, with three singular curves meeting at each such point.
    \item The idea of obtaining a relation between these two quotient varieties originates in the following observation. The group actions of $G_{81}$ on $\mP^5$ and $G_{27}$ on $\mP^2\times\mP^2\times\mP^1$ are toric and hence can be implemented using the software package PALP~\cite{Kreuzer:2002uu}, \cite{Braun:2012vh}. This allows one to easily verify that the polytopes of the two quotient ambient varieties have the same reflexive hull $P$ which admits a unique codimension two nef partition corresponding to a Calabi--Yau threefold with Hodge numbers $h^{11}=73, h^{21}=1$. This implies that there should be a birational equivalence between the families $\cX_{\lambda}$ and $\cY_z$. In the following, we will make this birational equivalence explicit.
\end{enumerate}
\end{remark}

In order to relate the varieties $\cX_\lambda$ and $\cY_z$, we introduce the following variety. 
\begin{definition}
Let $(u_0,u_1,u_2,v_0,v_1,v_2)$ be coordinates on $(\mC^\times)^6$. We define
\begin{equation}
 \mathcal T_w^\circ
:=\left\{\begin{matrix}
u_0+u_1+u_2=v_0+v_1+v_2=1,\quad
u_0u_1u_2v_0v_1v_2=w
\end{matrix}
\right\}\subset (\mC^\times)^6.
\end{equation}
\end{definition}
The relation between the varieties $\cW_z$ and $\cT_w^\circ$ is given by the following
\begin{lemma}
\label{prop:Y2T}
Let $\mathcal W_z^\circ\subset\mathcal W_z$ be the dense open subset on which $t\neq0,\infty$ and $p_xp_y\neq0$. 
Define
\[
u_i:=\frac{tx_i^3}{p_x},\qquad
v_j:=\frac{(z/t)y_j^3}{p_y},
\qquad 0\leq i,j\leq2.
\]
Then these functions define a morphism
\[
q_z^\circ:\mathcal W_z^\circ\longrightarrow
\mathcal T_{z^3}^\circ.
\]
The morphism $q_z^\circ$ is invariant under ${G_{27}}$, finite \'etale of degree $27$, and its geometric fibers are precisely the ${G_{27}}$-orbits. In particular,
\[
\mathcal W_z^\circ/{G_{27}}\cong \mathcal T_{z^3}^\circ.
\]
\end{lemma}

\begin{proof}
It follows immediately from the action~\eqref{eq:G27-action} that the functions $u_i,v_j$ are invariant. 
The defining equations of $\mathcal W_z^\circ$ give
\[
\sum_{i=0}^2u_i=\frac{ts_x}{p_x}=1,\qquad
\sum_{j=0}^2v_j=\frac{(z/t)s_y}{p_y}=1,
\]
and multiplication gives
\[
u_0u_1u_2=t^3,\qquad
v_0v_1v_2=\left(\frac zt\right)^3,
\]
so that $u_0u_1u_2v_0v_1v_2=z^3$. This proves that $q_z^\circ$ takes values in $\mathcal T_{z^3}^\circ$.

It remains to determine its fibers. Consider the homomorphism of algebraic tori
\begin{align*}
\vartheta:\{[w_0:w_1:w_2]\in\mathbb P^2:w_0w_1w_2\ne 0\}&\longrightarrow
\left\{(r_0,r_1,r_2)\in\mathbb A^3:r_0r_1r_2=1\right\},\\
[w_0:w_1:w_2]&\longmapsto\left(\frac{w_0^3}{w_0w_1w_2},
\frac{w_1^3}{w_0w_1w_2},
\frac{w_2^3}{w_0w_1w_2}\right)
\end{align*}
In the affine torus coordinates $a=w_1/w_0$, $b=w_2/w_0$, the last two coordinates of $\vartheta$ are $a^2b^{-1}$ and $a^{-1}b^2$. Thus the exponent matrix is
$\begin{pmatrix}2&-1\\-1&2\end{pmatrix}$, whose determinant is $3$. Therefore $\vartheta$ is a finite \'etale torus isogeny of degree $3$.

Now fix $(u_0,u_1,u_2,v_0,v_1,v_2)\in \mathcal T_{z^3}^\circ$. The equation
\[
t^3=u_0u_1u_2
\]
has exactly three solutions $t\in\mathbb C^\times$. For each such $t$, put $r_i=u_i/t$. Then $r_0r_1r_2=1$, so $(r_0,r_1,r_2)$ has exactly three inverse images $[x_0:x_1:x_2]$ under $\vartheta$. Similarly, since
\[
v_0v_1v_2=\frac{z^3}{u_0u_1u_2}=\left(\frac zt\right)^3,
\]
the triple $s_j=v_j/(z/t)$ satisfies $s_0s_1s_2=1$ and has exactly three inverse images $[y_0:y_1:y_2]$. The conditions $\sum u_i=\sum v_j=1$ then recover the two equations defining $\mathcal W_z^\circ$. Hence every geometric fiber of $q_z^\circ$ consists of exactly $3^3=27$ discrete points and $q_z^\circ$ has degree $27$.

Finally, the action of ${G_{27}}$ on $\mathcal W_z^\circ$ is free. Indeed, if an element fixes a point, then, since all $x_i$ and $y_j$ are nonzero, projective equality forces $a_0=a_1=a_2$ and $b_0=b_1=b_2$ in $\mathbb Z_3$. Such an element lies in the two projective scalar subgroups and hence represents the identity in ${G_{27}}$. Every ${G_{27}}$-orbit therefore has $27$ points. Since $q_z^\circ$ is ${G_{27}}$-invariant and every geometric fiber also has $27$ points, each fiber is exactly one ${G_{27}}$-orbit. Thus $q_z^\circ$ is the finite \'etale geometric quotient by ${G_{27}}$, and
\[
\mathcal W_z^\circ/{G_{27}}\cong \mathcal T_{z^3}^\circ
\]
as claimed. 
\end{proof}
In a similar manner, the relation between the varieties $\cV_\lambda$ and $\cT_w^\circ$ is given by the following
\begin{lemma}
\label{prop:X2T}
Suppose that $\lambda\ne 0,1$. 
Let $\mathcal V_\lambda^\circ\subset\mathcal V_\lambda$ be the dense open subset on which $p_xp_y\neq0$. 
Define
\[
u_i:=\frac{x_i^3}{3\lambda p_y},\qquad
v_j:=\frac{y_j^3}{3\lambda p_x},
\qquad 0\leq i,j\leq2.
\]
Then these functions define a morphism
\[
q_\lambda^\circ:\mathcal V_\lambda^\circ\longrightarrow
\mathcal T_{(3\lambda)^{-6}}^\circ
\]
The morphism $q_\lambda^\circ$ is invariant under ${G_{81}}$, finite \'etale of degree $81$, and its geometric fibers are precisely the ${G_{81}}$-orbits. In particular,
\[
\mathcal V_\lambda^\circ/{G_{81}}\cong \mathcal T_{(3\lambda)^{-6}}^\circ.
\]
\end{lemma}

\begin{proof}
The argument is parallel to the preceding lemma. The open subvariety $\mathcal V_\lambda^\circ$ is given by
\[
\frac{s_x}{3\lambda p_y}
=
\frac{s_y}{3\lambda p_x}
=1.
\]
It follows immediately from the action~\eqref{eq:G81-action} 
that the functions $u_i=x_i^3/(3\lambda p_y)$ and $v_j=y_j^3/(3\lambda p_x)$ are invariant. Moreover, the defining equations of $\mathcal V_\lambda^\circ$ give $\sum_i u_i=\sum_j v_j=1$, while multiplication gives $u_0u_1u_2v_0v_1v_2=(3\lambda)^{-6}$, so the resulting map takes values in $\mathcal T_{(3\lambda)^{-6}}^\circ$. As in the proof of Lemma~\ref{prop:Y2T}, this map is the restriction of a finite torus quotient: on the dense torus of $\mathbb P^5$, the corresponding monomial map has Smith normal form $\operatorname{diag}(1,1,3,3,9)$, hence it is a finite \'etale isogeny of degree $3\cdot3\cdot9=81$, with kernel $\mu_3^2\times\mu_9\cong\mathbb Z_3^2\times\mathbb Z_9$; the kernel is precisely the diagonal group $G_{81}$. Restricting this torus quotient to the equations $\sum_i u_i=\sum_j v_j=1$ is a base change, so it remains finite \'etale of degree $81$, and its geometric fibers are exactly the free $G_{81}$-orbits. Therefore $\mathcal V_\lambda^\circ/G_{81}\cong\mathcal T_{(3\lambda)^{-6}}^\circ$.
\end{proof}

\begin{proposition}
\label{prop:birational-quotients}
For nonzero parameters, the quotients have dense open subsets isomorphic to
$\cT_w^\circ$ with $w=(3\lambda)^{-6}$ for $\cX_\lambda=\cV_\lambda/G_{81}$ and $w=z^3$ for $\cY_z=\cW_z/G_{27}$.
In particular,
\begin{equation}\label{eq:birational-parameter-relation}
 \cX_\lambda\dashrightarrow \cY_z\quad\text{is birational if}\quad
 z^3=(3\lambda)^{-6}.
\end{equation}
\end{proposition}
\begin{proof}
This is a direct consequence of Lemmas~\ref{prop:Y2T} and~\ref{prop:X2T}.
\end{proof}

Indeed, for $z=1/(9\lambda^2)$ the correspondence could be written explicitly as
\begin{equation}
  \label{eq:V-to-W-rational-map}
 [x_0:x_1:x_2:y_0:y_1:y_2]\longmapsto
 \bigl([x_0:x_1:x_2],[y_0:y_1:y_2],[\px:3\lambda\py]\bigr).
\end{equation}
Substituting $t_0=\px$ and
$t_1=3\lambda\py$ into the defining equations
of $\mathcal W_z$ gives
\[
\begin{aligned}
t_0\sx-t_1\px
&=\px(\sx-3\lambda\py)=0,\\
zt_1\sy-t_0\py
&=(9\lambda^2z-1)\px\py=0,
\end{aligned}
\]
by \eqref{eq:Vlambda-equations} and
$z=1/(9\lambda^2)$.

Conversely, consider a point
$([x],[y],[t_0:t_1])$ in the dense open
subset $\mathcal W_z^\circ$ where
$t_0t_1\px\py\ne0$.
Choose representatives $x$ and $y$ for
the two projective blocks. Any inverse image
in $\mathcal V_\lambda$ must be represented
by
\[
[ax_0:ax_1:ax_2:by_0:by_1:by_2],
\qquad a,b\in\mathbb C^*.
\]
Since a common rescaling of $a$ and $b$
does not change the projective point, only
their ratio $s:=a/b$ matters. The defining
equations of $\mathcal V_\lambda$ then become
\[
s^3\sx=3\lambda\py,\qquad
\sy=3\lambda s^3\px.
\]
The equations of $\mathcal W_z$, together with
$z=1/(9\lambda^2)$, imply $\sx\sy=9\lambda^2\px\py$.
Hence the second equation follows from the
first. Thus there are exactly three inverse
choices, corresponding to the three roots of
$s^3\sx=3\lambda\py$. 

This map~\eqref{eq:V-to-W-rational-map} is equivariant for
\begin{equation}
  \label{eq:G81-to-G27-homomorphism}
  G_{81}\longrightarrow G_{27},\qquad
  (\alpha,\beta,\gamma)\longmapsto
  (\alpha-\gamma,-\beta+\gamma,\gamma\bmod3).
\end{equation}
This homomorphism is surjective with kernel
\begin{equation}
  \label{eq:K}
  K:=\ker(G_{81}\to G_{27}) \cong \{(0,0,3k):k\in\mZ_9\}\cong\mathbb Z_3\le G_{81}.     
\end{equation}
This kernel permutes the three block-scaling choices transitively. Thus it induces a degree-one map on the quotients, giving another proof of the birationality. 

For later purposes we record the following
\begin{lemma}
\label{prop:invariant-three-forms}
For smooth parameters $\lambda^6\ne 0,1$ and $z\not\in\Delta$, the actions of
$G_{81}$ on $\mathcal V_\lambda$ and $G_{27}$ on
$\mathcal W_z$ preserve their respective holomorphic
three-forms. 
\end{lemma}

\begin{proof}
For a smooth bicubic complete
intersection in $\mathbb P^5$ defined by $F_1=F_2=0$, suppose that the linear action by an element $g\in\mathrm{GL}_6(\mathbb C)$ fixes the equations $F_1,F_2$ up to scalar, respectively. Then on its Poincar\'e residue three-form $\Omega$ the action is given by
\[
g^*\Omega=\frac{\det(g)}
{\chi_{F_1}(g)\chi_{F_2}(g)}\,\Omega,
\]
where $\det(g)$ is the determinant of $g$ as a linear action and $\chi_{F_i}(g):=(g^*F_i)/F_i$ are the characters under which the defining equations transform. This follows from the transformation of the ambient volume form and the functoriality of the Poincaré residue.

For $\mathcal V_\lambda \subset \mP^5$, we have $F_1=s_x-3\lambda p_y,\ F_2=s_y-3\lambda p_x$ and the action \eqref{eq:G81-action} yields
\[
\begin{aligned}
\det(g)&=1,&
g^*F_1&=\xi_9^{3\gamma}F_1,&
g^*F_2&=\xi_9^{-3\gamma}F_2, 
\end{aligned}
\]
Thus the determinant character equals the product of the two equation characters, and the residue
three-form is $G_{81}$-invariant.

Similarly, for a smooth complete intersection $F_1=F_2=0$ in $\mathbb P^2\times\mathbb P^2\times\mathbb P^1$ of degree $(3,0,1)$ and $(0,3,1)$ suppose we are given $g\in\mathrm{GL}_3(\mathbb C)\times\mathrm{GL}_3(\mathbb C)\times\mathrm{GL}_2(\mathbb C)$ that fixes $F_1,F_2$ up to scalar. Then on its Poincar\'e residue three-form $\Omega$ the action is given by
\[
g^*\Omega=\frac{\det(g)}
{\chi_{F_1}(g)\chi_{F_2}(g)}\,\Omega
\]
again. So for $\mathcal W_z\subset
\mathbb P^2\times\mathbb P^2\times\mathbb P^1$ we have $F_1=t_0\sx - t_1\px,\ F_2=zt_1\sy - t_0\py$ and
the action \eqref{eq:G27-action} yields
\[
\begin{aligned}
\det(h)&=\xi_3^\gamma,&
h^*F_1&=\xi_3^\gamma F_1,&
h^*F_2&=F_2.    
\end{aligned}
\]
Therefore $G_{27}$ fixes the holomorphic three-form $\Omega_z$. 
\end{proof}

\subsection{Isomorphism of Hodge structures}
\label{sec:isom_Hodge}

The goal of this subsection is to study the map on the Hodge structure on the middle cohomology induced by this birational equivalence. For this purpose with start with two general lemmas on isomorphisms of Hodge structures under quotients and under certain blow--ups. 

\begin{lemma}
  \label{lem:finite-quotient-cohomology}
  Let a finite group $F$ act algebraically on a smooth
complex projective variety $T$, and let $q : T  \to Y=T/F$ be the
quotient morphism. Then
\begin{equation}
  \label{eq:finite-quotient-invariants}
  q^* : H^k(T/F) \xrightarrow{\sim} H^k(T)^F  
\end{equation}
is an isomorphism of rational Hodge structures.
If $F$ is normal in a larger finite group $G$ acting algebraically on
$T$, this isomorphism is $G/F$--equivariant.
\end{lemma}
\begin{proof}
  Set $\cA = q_*\mQ_T$. The map $q$ is finite, hence
  proper. Its fibers are discrete spaces. By proper base change, we
  have
  \[
    (R^j\cA)_y = H^j(q^{-1}(y),\mQ) = 0, \quad j>0.
  \]
  By the Leray spectral sequence, $H^k(Y,\cA)$ is naturally
  (and $F$--equivariantly) isomorphic to
  $H^k(T,\mQ)$.
  Since the stalk $\cA_y$ is a single orbit under $F$, its
  invariant functions are exactly the constants. Hence, $\mQ_Y
  \cong \cA^F$.

  Consider the averaging operator $e = \frac{1}{|F|}
  \sum_{f\in F} f$ acting on $\cA$. The decomposition $\mathcal
  A = \ker e \oplus \im e$ is preserved by taking cohomology. Hence, 
  taking invariants commutes with stalks and with taking cohomology. It
  follows that
  \[
    H^k(Y,\mQ) \cong H^k(Y,\cA)^F \cong H^k(T,\mQ)^F
  \]
  This composition is identified with the ordinary pullback $q^*$
  since it is induced by the unit $\mQ_Y \to q_*\mQ_T$
  followed by the Leray map.

  By~\cite[Theorem 5.33(iii)]{Peters:2008ab}, the map $q^*$ and every
  $f^*$ are morphisms of rational mixed Hodge structures. Therefore
  $H^k(T,\mQ)^F = \im e$ is a sub--mixed Hodge structure. The map
  $q^*$ onto $\im e$ is bijective. By~\cite[Corollary
  3.7]{Peters:2008ab}, it is therefore an isomorphism of mixed Hodge structures.

  If $F \triangleleft G$, then $q \circ g = \bar g \circ q$ gives
  $g^*q^* = q^*\bar g^*$ for every $g \in G$. This proves equivariance
  for $G/F$ as an isomorphism of mixed Hodge structures.
\end{proof}

\begin{remark}
Note that since $T$ is smooth projective, $H^k(T,\mQ)$ is pure
of weight $k$. Hence $H^k(Y,\mQ)$ is also pure of weight
$k$. This consequence does not require that the action of $F$ on $T$
be free.    
\end{remark}

\begin{lemma}
  \label{lem:blowup-H3-invariance}
  Let $b : \widetilde T \to T$ be the blow--up of a smooth complex
  projective threefold along a disjoint union of smooth centers, each
  a point or a projective line. Then
  \[
    b^* : H^3(T) \xrightarrow{\sim} H^3(\widetilde T)
  \]
  is an isomorphism of rational Hodge structures. It is equivariant
  for a finite group $F$ acting on $T$ preserving the union of the
  centers. 
\end{lemma}
\begin{proof}
  Let $C \subset T$ denote the smooth closed center in $T$ of codimension $r\geq 2$. Then we have the fiber square
  \[
    \begin{tikzcd}
        E=\mP(N_{C/T}) \ar[r,hook,"\iota"] \ar[d,"\pi"] &\widetilde T \ar[d,"b"] \\
        C \ar[r,hook,"j"] &T  \\
    \end{tikzcd}
  \]
  By the proper modification exact sequence of mixed Hodge
  structures~\cite[Corollary--Definition 5.37]{Peters:2008ab} and the fact that $T, \widetilde T, C$ and $E$ are all smooth projective, we obtain a short exact sequence of mixed Hodge structures
  \[
    0 \longrightarrow H^k(T) \xrightarrow{(b^*,j^*)} H^k(\widetilde T) \oplus H^k(C) \xrightarrow{\iota^*-\pi^*} H^k(E) \longrightarrow 0 \,.
  \]
  Next, put $\xi = c_1(\cO_E(1))$ using the convention that $\cO_E(-1)$ is the tautological line subbundle. Then the projective bundle formula is the isomorphism~\cite[Prop. 13.1]{Lewis:1999ab}
  \[
  \begin{aligned}
      \Phi_k : \bigoplus_{j=0}^{r-1} H^{k-2j}(C,\mQ)(-j) &\xrightarrow[]{\sim} H^k(E,\mQ)\\
      (\alpha_0, \dots, \alpha_{r-1}) & \mapsto \sum_{j=0}^{r-1} \xi^j \smile \pi^*\alpha_j
  \end{aligned}
  \]
  In particular, $b^*$ is injective. If $b^*\beta = 0$, then $\pi^*j^*\beta= \iota^*b^*\beta=0$, so $j^*\beta=0$. Since $(b^*,j^*)$ is injective, $\beta$ must be zero. Restriction therefore induces an isomorphism of Hodge structures
  \[
    \frac{H^k(\widetilde T)}{b^*H^k(T)} \xrightarrow[\iota^*]{\sim} \frac{H^k(E)}{\pi^*H^k(C)} \cong \bigoplus_{j=1}^{r-1} H^{k-2j}(C,\mQ)(-j)
  \]
  Note that the summand for $j=0$ is exactly the copy of $H^k(C)$ which is quotiented out. Next, we use the Gysin map and the self--intersection formula to write
  \[
    \iota^*\iota_*(\xi^{j-1} \smile \pi^*\alpha) = c_1(N_{E/\widetilde T}) \smile \xi^{j-1} \smile \pi^*\alpha = -\xi^j \smile \pi^*\alpha
  \]
  These classes therefore lift all the exceptional summands and give
  \[
  \begin{aligned}
    H^k(T) \oplus \bigoplus_{j=1}^{r-1} H^{k-2j}(C,\mQ)(-j) &\xrightarrow[]{\sim} H^k(\widetilde T)\\
    (\beta,\alpha_1,\dots,\alpha_{r-1}) & \mapsto b^*\beta - \sum_{j=1}^{r-1} \iota_*(\xi^{j-1} \smile \pi^*\alpha)
  \end{aligned}
  \]
  The exponent now is $j-1$ because $\iota_*$ has Hodge type $(1,1)$. Finally, setting $k=3$ and $C=\mathrm{point}$ or $C=\mP^1$ one sees that the extra cohomology groups vanish.
  Equivariance follows from naturality and the canonical lift of automorphisms preserving the center.
\end{proof}

Now, we return to the birational map $\cX_\lambda\dashrightarrow \cY_z$ in~\eqref{eq:birational-parameter-relation}. In order to make use of these lemmas, we decompose this map into a sequence of simpler quotient maps and blow--up maps. It will turn out that these simpler maps are naturally given by the geometry and the group actions that we have chosen. 

In the following until the end of this section we simplify the notation and write $V=\cV_\lambda$, $W=\cW_z$, $X=\cX_\lambda$ and $Y=\cY_z$. Also, all the cohomology groups will be singular cohomology groups even when the topological space in question is singular. 

To start we set
\[
  \begin{aligned}
    B_x &= \{\sx = p_x = 0\} \subset \mP^2_x, & B_y = \{\sy =
    p_y = 0\} \subset \mP^2_y.
  \end{aligned}
\]
Each set consists of nine reduced points. At any such point exactly one coordinate vanishes.
The product has a nonzero derivative in that coordinate direction, while the cubic sum has a
nonzero derivative in the independent ratio coordinate. Hence
$(\sx,p_x)$, respectively $(\sy,p_y)$, are local analytic coordinates
at these points.
\begin{proposition}
  \label{prop:K-quotient-H3}
  Set $P = \{ (x:0) \mid x \in B_x \} \sqcup \{ (0:y) \mid y \in B_y
  \}$ and $b: \widetilde V := \Bl_PV\to V$.  
  The eighteen points $P$ are exactly the fixed locus of every nontrivial element of $K$ defined in~\eqref{eq:K}. The quotient
  \[
    q: \widetilde V \to Z := \widetilde V/K
  \]
  is finite of degree three, and $Z$ is smooth and projective. The lifted $G_{81}$--action induces an
$G_{27} \cong G_{81}/K$--action on $Z$, with a natural $G_{27}$--equivariant identification of rational Hodge structures 
\[
  H^3(Z) \simeq H^3(V)^K\ .
\]
\end{proposition}
\begin{proof}
  A generator of $K$ acts as $(x:y) \mapsto (\xi_3 x:
  \xi_3^{-1}y)$. The fixed locus consists of $\mP^2_x \times \point$
  and $\point \times \mP^2_y$. The intersection with $V$ gives the set $P$. At $(x: 0)$, the
independence of the derivatives of $s_x$ and $p_x$ solves the two
local $x$ coordinates in terms of the
three $y$ coordinates. Relative to a nonzero $x$ coordinate, $K$ acts
on these analytic coordinates as multiplication by
$\xi_3$. Similarly, $K$ acts on the analytic coordinates at $(0:y)$
by $\xi_3^{-1}$. In a blowup chart $(s,su,sv)$ the action is $(s,u,v)
\mapsto (\xi_3^{\pm 1}s,u,v)$, whose invariant coordinate ring is
$\mathbb C[s^3,u,v]$. Thus the quotient is smooth along every
exceptional divisor; elsewhere the action is free. Projectivity follows from taking a finite quotient of a
projective variety. The set $P$ is $G_{81}$--invariant, so the actions lift
and descend as asserted. Apply Lemma~\ref{lem:blowup-H3-invariance} to $b$ and
Lemma~\ref{lem:finite-quotient-cohomology} to $q$ to obtain the final statement.
\end{proof}

\begin{proposition}
For $\lambda^6\ne 0,1$, consider the following hypersurface of bidegree $(3,3)$ in $\mP^2_x \times \mP^2_y$,
\[
    B = \{\sx\sy = 9\lambda^2p_xp_y\} \subset \mP^2_x \times \mP^2_y\ .
\]
Then $B$ is a singular Calabi--Yau threefold whose singular locus is $N = B_x \times B_y$, and every point of $N$ is a conifold singularity of $B$. Assume that 
\begin{equation}
  \label{eq:small-resolution-parameters}
  z=(9\lambda^2)^{-1}\ .
\end{equation}
Then there are $G_{27}$--equivariant projective birational morphisms
\[
  \begin{aligned}
    \rho_Z : Z &\to B, & \rho_W : W \to B
  \end{aligned}
\]
which are isomorphisms away from $N$. Over every node they are the two opposite small
resolutions, with exceptional fiber $\mP^1$ of normal bundle
$\cO(-1) \oplus  \cO(-1)$.    
\end{proposition}
\begin{proof}
  The idea is simply to notice that the equation for $B$ is the vanishing of a determinant of a $2\times 2$ matrix so there will be conifold singularities at the zeros of this determinant. We then only need to show that the two simultaneous small resolutions precisely yield the varieties $Z$ and $W$. The morphism $\rho_W: W \to B$ is already discussed in~\cite{sch88}.

  First, $\rho_W$ forgets $t$. The two equations for $[t_0 : t_1]$ are the kernel equations for the
  matrix
  \[
    \begin{pmatrix}
      \sx & -p_x \\
      -p_y & z\sy
    \end{pmatrix}.
  \]
Its determinant vanishes on $B$; outside $N$ its rank is one. Its
kernel is then a unique projective point varying regularly on the
rank-one charts. Thus $W \to B$ is an isomorphism outside $N$ and has
fiber $\mP^1$ over $N$. Since $W$ is smooth, $B$ is smooth outside $N$
as well. At a point of $N$ use analytic coordinates
\[
  \begin{aligned}
    u&= \sx, & v &= p_x, & a &= \sy, & w &= 9\lambda^2p_y.    
  \end{aligned}
\]
The local equation is $ua= vw$, a conifold singularity, and $W$ has
incidence relations
\begin{equation}
    \label{eq:W-small-resolution}
  \begin{aligned}
    t_0u &= t_1v, & t_0w &= t_1a.
  \end{aligned}    
\end{equation}
On $\mP^2\times\mP^2$ set $L_x = \cO_{\mathbb
  P^2\times \mP^2 }(-1,0)$ and $L_y = \cO_{\mP^2\times
  \mP^2 }(0,-1)$ and consider the projective bundle $P_1=\mP_{\mP^2\times
  \mP^2 }(L_x \oplus L_y)$. It comes with two maps, the projection $p:P_1 \to \mP^2\times \mP^2$ and the map $q:P_1 \to \mP^5$ given by $([x],[y],[r:s]) \mapsto [rx:sy]$ in local trivializations. The map $q$ is a resolution of the rational map $\rho: \mP^5 \to \mP^2\times \mP^2$. Viewing $P_1$ as the incidence variety 
  \[
    \{([X:Y],([x],[y])) \mid X \in \mC x, Y \in \mC y\} \subset \mP^5 \times \mP^2 \times \mP^2,
  \]
  it is the closure of the graph of $\rho$.
  Its restriction to the graph closure $\rho|_{V\setminus P}$ is $\widetilde V$: near $[\alpha:0]$, the three $y$--coordinates are local coordinates on $V$, and the closure of the graph of $[y_0:y_1:y_2]$ is their blow--up. Near $[0:\beta]$, the same argument uses the three $x$--coordinates. Away from $P$, both blocks are nonzero and the graph needs no modification.
  The strict transform in $P_1$ has equations
\begin{equation}
  \label{eq:V-incidence-equations}
  r^3\sx=3\lambda s^3\py,\qquad s^3\sy=3\lambda r^3\px.
\end{equation}
These glue under rescaling $x,y$, since $r,s$ transform by inverse scalars.
At $r=0$ they give $y\in B_y$ with $x$ arbitrary, and at $s=0$ they give
$x\in B_x$ with $y$ arbitrary: precisely the eighteen exceptional planes.
Near each of these loci the two independent base functions can be solved
in terms of $r/s$ or $s/r$ and the other block coordinates. Thus the equations
introduce no extra component or thickening there.

A generator of $K$ acts as $(x:y) \mapsto (\xi_3 x:
  \xi_3^{-1}y)$. It induces a fiberwise $K$-action on $P_1$ by $[r:s]\mapsto[\xi_3 r:\xi_3^{-1}s]$.
Its quotient is the map to $P_3=\mP_{\mP^2\times \mP^2}(L_x^{\otimes3}\oplus L_y^{\otimes3})$,
with $[R:S]=[r^3:s^3]$. The transition functions are cubed, so this description
is global. Taking invariants in \eqref{eq:V-incidence-equations} gives
\begin{equation}
  \label{eq:Z-incidence-model}
  Z=\{R\sx=3\lambda S\py,\quad S\sy=3\lambda R\px\}
  \subset P_3.
\end{equation}
Indeed, on $s\ne0$, the polynomial algebra in $r/s$ is free over the algebra
in $(r/s)^3$, with basis $1,r/s,(r/s)^2$. Both defining relations lie in
that invariant algebra; hence their quotient has invariant ring exactly as
in \eqref{eq:Z-incidence-model}. The other chart is identical.
Projection to $\mP_x^2\times \mP_y^2$ now gives $\rho_Z$. Its two linear equations in $[R:S]$
have rank one outside $N$ and rank zero at $N$, proving the same fiber
assertions as for $W$.

Near a node put $[l_0:l_1]=[S:3\lambda R]$. Then \eqref{eq:Z-incidence-model}
reads
\begin{equation}\label{eq:Z-small-resolution}
 l_0w=l_1u,\qquad l_0a=l_1v.
\end{equation}
For the rank-one matrix $\left(\begin{smallmatrix}u&v\\w&a\end{smallmatrix}\right)$,
\eqref{eq:W-small-resolution} records the common row direction and
\eqref{eq:Z-small-resolution} the common column direction. These are its two opposite
small resolutions. Their affine charts identify neighborhoods of the
exceptional curve with neighborhoods of the zero section in
$\cO_{\mP^1}(-1)^{\oplus2}$.
Finally the block action of $G_{81}$ factors through \eqref{eq:G81-to-G27-homomorphism}; on $B$
the cubic sums are invariant and the products have inverse characters.
Consequently both morphisms are $G_{27}$-equivariant.
\end{proof}

\begin{proposition}
\label{prop:common-blowup}
There is a smooth projective $G_{27}$-variety
\[
 M=\operatorname{Bl}_N B
\]
with $G_{27}$-equivariant morphisms $\sigma_Z:M\to Z$ and $\sigma_W:M\to W$.
Each is the blowup of the disjoint union of its $81$ exceptional projective
lines over $N$. In particular,
\begin{equation}
  \label{eq:common-blowup-H3}
  \begin{aligned}
     \sigma_Z^*:H^3(Z) &\xrightarrow{\sim}H^3(M), & \sigma_W^*:H^3(W) &\xrightarrow{\sim}H^3(M)
  \end{aligned}
\end{equation}
is a pair of $G_{27}$-equivariant isomorphisms of rational Hodge structures.
\end{proposition}
\begin{proof}
The assertion is local over $N$. For the node $ua=vw$, record both the row
and column directions of its rank-one matrix. The resulting space is
$\operatorname{Tot}(\cO_{\mP^1\times\mP^1}(-1,-1))$.
It is the blowup of the vertex: the projectivized rank-one cone is the
quadric $\mP^1\times\mP^1\subset\mP^3$, and its tautological line bundle is
$\cO_{\mP^1\times\mP^1}(-1,-1)$. Forgetting one direction contracts one ruling of the
quadric and is the blowup of the zero section in the corresponding
$\cO_{\mP^1}(-1)^{\oplus2}$ small resolution.
For example, on $t_0\ne0$ in \eqref{eq:W-small-resolution}, with $t=t_1/t_0$,
the pullback of the node ideal $(u,v,w,a)$ is $(v,a)$, the ideal of the
exceptional curve. The other charts give the same conclusion.
The universal property of blowup glues the local identifications; away from
$N$ every map is the identity on $B$. The resulting $M$ is smooth and
projective. All centers are $G_{27}$-stable, making these identifications
equivariant. Lemma~\ref{lem:blowup-H3-invariance} proves \eqref{eq:common-blowup-H3}.
\end{proof}


\begin{theorem}\label{thm:quotient-H3-comparison}
Assume \eqref{eq:small-resolution-parameters}. Set $R=M/G_{27}$. There are proper birational
morphisms
\begin{equation}
  \label{eq:birational-roof}
  \cX_\lambda\xleftarrow{\ f\ }R\xrightarrow{\ g\ } \cY_z
\end{equation}
whose pullbacks on the third singular cohomology group are isomorphisms. Thus
\begin{equation}\label{eq:quotient-H3-isomorphism}
 \Phi_\lambda=(g^*)^{-1}f^*:
 H^3(\cX_\lambda,\mathbb Q)\xrightarrow{\sim}H^3(\cY_z,\mathbb Q)\ 
\end{equation}
is an isomorphism of pure rational Hodge structures of weight three.
\end{theorem}
\begin{proof}
The map $b:\widetilde V\to V$ descends to $Z\to V/K$; composing with
$V/K\to V/G_{81}=X$ gives an $G_{27}$-invariant morphism $j:Z\to X$.
Thus $j\circ \sigma_Z$ descends to $f:R\to X$ and the composite $M\xrightarrow{\sigma_W} W\to Y$
descends to $g:R\to Y$. Both are proper and birational. On the dense open subset
where these modifications are isomorphisms, their correspondence is exactly
\eqref{eq:V-to-W-rational-map}, so they realize the birational matching of
\cref{prop:birational-quotients}.
Write $q_X:V\to X$, $q_Y:W\to Y$, and $q_M:M\to R$ for the quotient maps.
We can visualize all maps in the following commutative diagram of $G_{27}$--equivariant morphisms:
\[
  \begin{tikzcd}
     & & \widetilde V \ar[r,"b"] \ar[d,"q"] & V \ar[d, "q_X"] \\
     & M \ar[r,"\sigma_Z"] \ar[d,"\sigma_W" ']  & Z=\widetilde V/K \ar[r,"j"] \ar[d,"\rho_Z"] &  X = V/G_{81} \\
     & W \ar[r,"\rho_W", near end ] \ar[dd,"q_Y" '] & B & & \\
     & & R=M/G_{27} \ar[dl,"g"] \ar[uur,"f" '] \ar[from=uul,"q_M" ', near end, crossing over]& & \\
     & Y=W/G_{27} & & \\
  \end{tikzcd}
\]
Since $j\circ q=q_X\circ b$, on cohomology
\[
 q^*j^*=b^*q_X^*.
\]
Lemma~\ref{lem:finite-quotient-cohomology} identifies $H^3(X)$ with $H^3(V)^{G_{81}}$ and identifies
$H^3(Z)^{G_{27}}$ with $(H^3(\widetilde V)^K)^{G_{27}}=H^3(\widetilde V)^{G_{81}}$.
By Lemma~\ref{lem:blowup-H3-invariance}, the point blowup $b$ induces an isomorphism between these $G_{81}$-invariant spaces. Therefore $j^*:H^3(X)\to H^3(Z)^{G_{27}}$ is an isomorphism.
Now
\[
 q_M^*f^*=\sigma_Z^*j^*,\qquad q_M^*g^*=\sigma_W^*q_Y^*.
\]
Lemma~\ref{lem:finite-quotient-cohomology} and \eqref{eq:common-blowup-H3} show that both right-hand sides
are isomorphisms onto $H^3(M)^{G_{27}}$. Hence $f^*$ and $g^*$ are isomorphisms
on \emph{ordinary} third cohomology, proving \eqref{eq:quotient-H3-isomorphism}.
Equivalently the same map is
\begin{equation}\label{eq:cover-H3-comparison}
 \Phi_\lambda=(q_Y^*)^{-1}(\sigma_W^*)^{-1}\sigma_Z^*(q^*)^{-1}b^*q_X^*,
\end{equation}
where each inverse quotient pullback is taken only on the invariant image
identified above.
The smooth projective covers $V,W,M$ imply purity of $H^3(X),H^3(Y),H^3(R)$.
All pullbacks are morphisms of mixed Hodge structures; their bijections have
Hodge-theoretic inverses by strictness. This proves the Hodge assertion.
\end{proof}

In the rest of this section, put $\lambda_*=\sqrt{-2}$, $z_*=-1/18$ and $w_*=z_*^3=(3\lambda_*)^{-6}$. 
Note that $W:=\mathcal W_{z_*}, Y:=\mathcal Y_{z_*}$ are defined over $\mathbb Q$, however, 
$\lambda_*=\sqrt{-2}\not\in\mathbb Q$. 
To prove the arithmetic statements for the threefolds $V:=\mathcal V_{\lambda_*},X:=\mathcal X_{\lambda_*}$, we need to modify the equations of $V, X$ to define them over $\mathbb Q$. 
These $\mathbb Q$-models $V^{(0)}, X^{(0)}$ will turn out to be useful in \cref{sec:point-counting-modularity} to investigate modularity. 

\begin{lemma}
\label{lem:rational-resolved-comparison}
The quotient $X=\mathcal X_{\lambda_*}$ has the following $\mathbb Q$-model:
\[
 X^{(0)}:=V^{(0)}/G_{81},\qquad
 V^{(0)}:=\{2\sr+3p_y=\sy-12p_r=0\}\subset\mathbb P^5_{\mathbb Q},
\]
where $\sr=\sum r_i^3$ and $p_r=\prod r_i$.
The birational roof of \cref{thm:quotient-H3-comparison}
descends to $\mathbb Q$:
\[
 X^{(0)}\xleftarrow{\ f\ }R=M/G_{27}
 \xrightarrow{\ g\ }Y.
\]
Its pullbacks on $H^3$ are isomorphisms, in both Betti and
$\ell$-adic cohomology.
\end{lemma}

\begin{proof}
The substitution $x_i=\lambda_*r_i$ identifies
$V^{(0)}_{\mathbb Q(\sqrt{-2})}$ with
$V=V_{\lambda_*}$, compatible with the $G_{81}$-action. 
By \cref{prop:X2T}, the invariant torus coordinates on $V^{(0)}$ are
\[
 u_i=-\frac{2r_i^3}{3p_y},\qquad
 v_j=\frac{y_j^3}{12p_r}.
\]
They identify a dense open subset of $X^{(0)}$ over $\mathbb Q$
with $\mathcal T_{w_*}^{\circ}$, as \cref{prop:Y2T} does for $Y$. 
The common blowup $M$ and $R$
are defined over $\mathbb Q$; the morphism
$R_{\mathbb Q(\sqrt{-2})}\to X$ of
\cref{thm:quotient-H3-comparison} descends to $R\to X^{(0)}$,
since its Galois conjugates agree on this common dense open
and the target is separated.
Finally the Betti assertion is \cref{thm:quotient-H3-comparison};
functorial Betti--\'etale comparison yields the $\ell$-adic one. 
\end{proof}

\begin{theorem}[Smooth projective crepant resolutions over $\mathbb Q$]
\label{thm:Q-crepant-resolutions}
Let $X^{(0)}=V^{(0)}/G_{81}$ and $Y=W/G_{27}$ be the
$\mathbb Q$-models of \cref{lem:rational-resolved-comparison}.
There exist smooth projective crepant resolutions over $\mathbb Q$,
\[
 \rho_X:\widetilde X\longrightarrow X^{(0)},\qquad
 \rho_Y:\widetilde Y\longrightarrow Y.
\]
Their sources are smooth projective Calabi--Yau threefolds.
After base change to $\mathbb Q(\lambda_*)$, the first morphism
resolves $\mathcal X_{\lambda_*}$; the second resolves
$\mathcal Y_{z_*}$. 
\end{theorem}

\begin{proof}
The diagonal actions \eqref{eq:G81-action} and
\eqref{eq:G27-action} are given by monomial characters of roots
of unity. Hence they descend to actions of finite diagonalizable
$\mathbb Q$-group schemes
\[
 G_{81}\simeq\mu_3^2\times\mu_9,\qquad
 G_{27}\simeq\mu_3^3,
\]
where $\mu_n$ denotes the group scheme of $n$th roots of unity
over $\mathbb Q$. These group schemes act on $V^{(0)}$ and $W$,
respectively, with the indicated quotients.
Both covers are smooth projective threefolds:
for $V^{(0)}$ this follows from its identification with
$\mathcal V_{\lambda_*}$ after extending scalars to
$\mathbb Q(\lambda_*)$, and for $W$ from $(9z_*)^3\ne1$.

By \cref{prop:invariant-three-forms}, 
the two actions preserve their nowhere-vanishing holomorphic
three-forms. So the canonical bundle
of each cover is trivial as an equivariant line bundle.
By the construction of $G$-Hilbert schemes for finite group
schemes \cite{Blume2011}, the $\mathbb Q$-schemes
\[
\begin{aligned}
 \widetilde X:=G_{81}\text{-}\operatorname{Hilb}\bigl(V^{(0)}\bigr),\qquad
 \widetilde Y:=G_{27}\text{-}\operatorname{Hilb}(W)
\end{aligned}
\]
exist and carry projective Hilbert--Chow morphisms to the
respective quotients. Their formation is compatible with
extension of the ground field. 

Since the covers are smooth Calabi--Yau threefolds with equivariantly trivial canonical bundles, the theorem of Bridgeland--King--Reid \cite[Theorem~1.2]{BKR2001} shows that both Hilbert--Chow morphisms are smooth crepant resolutions. Smoothness and crepancy descend to $\mathbb Q$, and the sources are projective since their targets are projective. 
\end{proof}

\begin{remark}
\begin{enumerate}
    \item In Remark~\ref{rem:toric} we observed that the two polytopes of $\mP^5/G_{81}$ and $\mP^2\times\mP^2\times\mP^1/G_{27}$ have the same reflexive hull $P$. Using standard methods of toric geometry, it is possible to construct maximal projective crepant partial desingularizations $\Sigma_{\widetilde X}$ and $\Sigma_{\widetilde Y}$ of the face fan $\Sigma_P$ of $P$, yielding explicit realizations of the smooth varieties $\widetilde \cX_{\lambda}$ and $\widetilde \cY_z$.
    \item We end with an overview of the Hodge numbers of the smooth Calabi--Yau threefolds constructed in this section.
\[
\begin{array}{@{}|l|c|c|r|r|@{}}
\hline & h^{1,1} & h^{2,1} \\\hline
V=\mathcal V_\lambda
  & 1 & 73 \\
W=\mathcal W_z
  & 19 & 19 \\
\widetilde X = \widetilde \cX_\lambda
  & 73 & 1 \\
\widetilde Y = \widetilde \cY_z
  & 73 & 1 \\
Z=\widetilde V/K
  & 19 & 19 \\\hline
\end{array}
\]
\end{enumerate}
\end{remark}

\section{Group actions on the middle cohomology}
\label{sec:special-fiber-symmetries}
\label{sec:middle-cohomology-actions}

In this section, we mainly specialize to the smooth fiber
\[
W:=\mathcal W_{z_*},\qquad z_*:=-\frac1{18}.
\]
This corresponds to the conjectured attractor point since
\[
w_* = z_*^3=-\frac{1}{5832}. 
\]
Since $z_*\notin\Delta$, the threefold $W$ is smooth and projective, with
$h^{1,1}(W)=h^{2,1}(W)=19$ and $\dim H^3(W,\mathbb C)=40$.  

We will firstly study the automorphism subgroups $\Gamma\cong\operatorname{SL}(2,\mathbb F_3)$ and $G_{27}$ of the special fiber $W=\mathcal W_{z_*}$ and compare their actions on its middle cohomology. Here $\Gamma$ comes from automorphisms of the Hesse pencil. 
This yields a decomposition of $G_{27}$-invariant subspace into two rank-two rational Hodge structures, which we show is motivic over $\mathbb Q$. 
We then transfer this decomposition to the crepant resolutions of the quotient families and establish their rank-two attractor property. 

\subsection{The special automorphism}
\label{sec:special_aut}

For this special fiber $W$, we can define an extra automorphism called $\psi$. 
After this, we will study the automorphism subgroup $\Gamma$ generated by $\psi$ and a diagonal action $\phi$ and their actions on the middle cohomology $H^3$. 

\begin{proposition}
\label{prop:special-fiber-automorphisms}
The following maps define automorphisms of $W$:
\begin{align}
\nonumber
\phi:\quad
&((x_0:x_1:x_2),(y_0:y_1:y_2),t)\\
\label{eq:special-phi-action}
&\longmapsto
((\xi_3x_0:x_1:x_2),(\xi_3^2y_0:y_1:y_2),\xi_3t),
\\[4pt]
\nonumber
\psi:\quad
&((x_0:x_1:x_2),(y_0:y_1:y_2),t)\\
\nonumber
&\longmapsto
\bigl((x_0+x_1+x_2:x_0+\xi_3x_1+\xi_3^2x_2:
 x_0+\xi_3^2x_1+\xi_3x_2),\\
\label{eq:special-psi-action}
&\hspace{2.2cm}
(y_0+y_1+y_2:y_0+\xi_3^2y_1+\xi_3y_2:
 y_0+\xi_3y_1+\xi_3^2y_2),f(t)\bigr).
\end{align}
Here
\begin{equation}
\label{eq:special-psi-base-map}
f(t):=\frac{1-3t}{3+18t}
\end{equation}
is an involution of $\mathbb P^1_t$. The subgroup
\[
\Gamma:=\langle\phi,\psi\rangle\subset\operatorname{Aut}(W)
\]
is isomorphic to $\operatorname{SL}(2,\mathbb F_3)$.
\end{proposition}

\begin{proof}
On the first Hesse factor the two maps are the automorphisms $\phi,\psi$ of
\cref{prop:hesse-automorphisms}, while on the second factor they are
$\phi^{-1},\psi^{-1}$.  For $\phi$, the two Hesse parameters are multiplied
by $\xi_3$ and $\xi_3^2$, respectively, so their product remains $z_*$.  For
$\psi$, a direct calculation gives
\[
f(t)f(z_*/t)=z_*.
\]
Thus both maps preserve the fiber-product equation.  
They satisfy the same relations as in \eqref{eq:gamma-relations}, namely
\[
\phi^3=\psi^4=(\phi\psi)^3=[\phi,\psi^2]=1. 
\]
Projecting to the first Hesse factor implies that the action is faithful. 
Hence they generate exactly the binary tetrahedral group $\operatorname{SL}(2,\mathbb F_3)$. 
\end{proof}

\begin{remark}
The automorphism $\phi$ exists on every fiber $\mathcal W_z$, as first discussed in Section 5.1 of \cite{bon23}. 
The second automorphism is special to $z=z_*$.  Indeed,
\[
f(t)f(z/t)-z
=-\frac{(18z+1)(3t^2+18tz-t+3z)}{9(t+6z)(6t+1)},
\]
so $f(t)f(z/t)=z$ identically in $t$ if and only if $z=-1/18$. 
\end{remark}

To study the representations, 
we use the following character table for $\operatorname{SL}(2,\mathbb F_3)$.  
The seven conjugacy classes of $\Gamma$ are denoted by
$\mathsf 1,\mathsf 2,\mathsf{3A},\mathsf{3B},\mathsf 4,\mathsf{6A},\mathsf{6B}$.
The elements $\phi,\phi^2,\phi\psi^2,\phi^2\psi^2$ belong to the classes
$\mathsf{3A},\mathsf{3B},\mathsf{6A},\mathsf{6B}$, respectively;
$\psi^2$ belongs to $\mathsf 2$, and $\psi$ belongs to $\mathsf 4$.  With
$\xi_6:=-\xi_3^2=(1+\sqrt{-3})/2$, the irreducible characters are
\begin{equation}
\label{eq:SL2F3-character-table}
\begin{array}{|c|ccccccc|}
\hline
\text{class} & \mathsf 1 & \mathsf 2 & \mathsf{3A} & \mathsf{3B}
& \mathsf 4 & \mathsf{6A} & \mathsf{6B} \\
\hline
\text{size} & 1 & 1 & 4 & 4 & 6 & 4 & 4 \\
\hline
\rho_1 & 1 & 1 & 1 & 1 & 1 & 1 & 1 \\
\rho_2 & 1 & 1 & \xi_3^2 & \xi_3 & 1 & \xi_3^2 & \xi_3 \\
\rho_3 & 1 & 1 & \xi_3 & \xi_3^2 & 1 & \xi_3 & \xi_3^2 \\
\rho_4 & 2 & -2 & -1 & -1 & 0 & 1 & 1 \\
\rho_5 & 2 & -2 & \xi_6^5 & \xi_6 & 0 & \xi_3 & \xi_3^2 \\
\rho_6 & 2 & -2 & \xi_6 & \xi_6^5 & 0 & \xi_3^2 & \xi_3 \\
\rho_7 & 3 & 3 & 0 & 0 & -1 & 0 & 0 \\
\hline
\end{array}
\end{equation}

\begin{lemma}
\label{lem:gamma-fixed-loci}
For $g\in \Gamma$, let
\[
\Lambda(g):=\sum_{k=0}^6(-1)^k
\operatorname{Tr}\bigl(g\mid H^k(W,\mathbb C)\bigr).
\]
Then
\[
\begin{array}{c|ccccccc}
\text{class} & \mathsf 1 & \mathsf 2 & \mathsf{3A} & \mathsf{3B}
& \mathsf 4 & \mathsf{6A} & \mathsf{6B} \\
\hline
\Lambda & 0 & 0 & 18 & 18 & 8 & 6 & 6.
\end{array}
\]
\end{lemma}

\begin{proof}
We compute the fixed loci for $\psi,\psi^2,\phi$ and $\phi\psi^2$.

The fixed points of the base involution $f$ are
\[
\tau_\pm:=\frac{-1\pm\sqrt3}{6},\qquad
\frac{z_*}{\tau_\pm}=\tau_\mp.
\]
At $t=\tau_\pm$, the fixed points on the first factor are
$(1\pm\sqrt3:1:1)$ and $O$, while those on the second factor are
$(1\mp\sqrt3:1:1)$ and $O$.  Hence $\psi$ has eight isolated fixed points,
so $\Lambda(\psi)=8$.

The fixed fibers of $\phi\psi^2$ are $t=0,\infty$.  Over $t=\infty$ the
fixed points are
\[
(x_0:x_1:x_2)=O,
\qquad
(y_0:y_1:y_2)\in\{(1:0:0),(0:1:1),(0:1:-1)\},
\]
and over $t=0$ the roles of the two factors are reversed.  Thus
$\phi\psi^2$ has six isolated fixed points.  The same holds for its inverse
$\phi^2\psi^2$, and therefore
$\Lambda(\phi\psi^2)=\Lambda(\phi^2\psi^2)=6$.

For $\phi$, the fixed fibers are again $0$ and $\infty$.  Over $t=\infty$
the first factor is the Fermat cubic
$x_0^3+x_1^3+x_2^3=0$, and the fixed points of
$x_0\mapsto\xi_3x_0$ are
\[
(0:-1:\xi_3^i),\qquad i=0,1,2.
\]
The second factor is the $I_3$-fiber $y_0y_1y_2=0$, on which the fixed locus
of $y_0\mapsto\xi_3^2y_0$ is the line $y_0=0$ together with the point
$(1:0:0)$.  Hence the contribution over $t=\infty$ is
$3(\chi(\mathbb P^1)+1)=9$.  The same calculation over $t=0$, with the two
factors interchanged, gives another $9$.  Thus
$\Lambda(\phi)=18$, and the same value holds for $\phi^2$.

Finally, $\psi^2$ acts on every Hesse fiber by
$(x_0:x_1:x_2)\mapsto(x_0:x_2:x_1)$, and similarly on the second factor.
Its fixed locus on $W$ consists of the following four curves:
\[
\begin{cases}
(x_0:x_1:x_2)=O,\quad (y_0:y_1:y_2)=O,\\[2pt]
(x_0:x_1:x_2)=O,\quad (y_0:y_1:y_2)=(y:1:1),\quad
\dfrac{z_*}{t}=\dfrac{y}{y^3+2},\\[7pt]
(x_0:x_1:x_2)=(x:1:1),\quad (y_0:y_1:y_2)=O,\quad
 t=\dfrac{x}{x^3+2},\\[7pt]
(x_0:x_1:x_2)=(x:1:1),\quad (y_0:y_1:y_2)=(y:1:1),\quad
 t=\dfrac{x}{x^3+2},\quad
\dfrac{z_*}{t}=\dfrac{y}{y^3+2}.
\end{cases}
\]
The first three curves are rational.  The fourth one is the compactification of
\begin{equation}
\label{eq:psi-squared-fixed-curve}
\mathcal C:\quad
\frac{x}{x^3+2}\cdot\frac{y}{y^3+2}=-\frac1{18},
\end{equation}
or equivalently
$(x^3+2)(y^3+2)+18xy=0$.  Its closure in
$\mathbb P^1\times\mathbb P^1$ has bidegree $(3,3)$.  A direct Jacobian check
in the four standard affine charts shows that it is smooth, hence it has genus $4$.  Therefore
\[
\Lambda(\psi^2)=3\chi(\mathbb P^1)+\chi(\mathcal C)
=6+(2-8)=0.
\]
Since $\chi(W)=0$, also $\Lambda(1)=0$.
\end{proof}

Let $F$ be the class of a fiber of $W\to\mathbb P^1$, and let
\[
D_x:=\bigl[\{O\}\times_{\mathbb P^1}\mathcal E_{z_*/t}\bigr],
\qquad
D_y:=\bigl[\mathcal E_t\times_{\mathbb P^1}\{O\}\bigr].
\]
These three classes are fixed by $\Gamma$.  We now choose the remaining component
classes in a way compatible with the present Hesse parameter.

\begin{definition}
\label{def:singular-fiber-component-classes}
Put
\[
\Sigma_x:=\left\{0,\frac13,\frac{\xi_3}{3},\frac{\xi_3^2}{3}\right\},
\qquad
\Sigma_y:=\left\{\infty,-\frac16,-\frac{\xi_3}{6},-\frac{\xi_3^2}{6}\right\}.
\]
For $u\in \Sigma_x\setminus\{0\}$ set $\lambda_u=(3u)^{-1}$.  The three
components of the singular fiber $\mathcal E_u$ are
\begin{align*}
L_u^O&:\ \lambda_ux_0+x_1+x_2=0,\\
L_u^+&:\ \lambda_ux_0+\xi_3x_1+\xi_3^2x_2=0,\\
L_u^-&:\ \lambda_ux_0+\xi_3^2x_1+\xi_3x_2=0.
\end{align*}
At $u=0$ we put
\[
L_0^O=\{x_0=0\},\qquad L_0^+=\{x_1=0\},\qquad L_0^-=\{x_2=0\}.
\]

For $v\in \Sigma_y\setminus\{\infty\}$ set $\mu_v=-6v$.  Since the second
Hesse parameter is $z_*/v$, the three components of its singular fiber are
\begin{align*}
M_v^O&:\ \mu_vy_0+y_1+y_2=0,\\
M_v^+&:\ \mu_vy_0+\xi_3y_1+\xi_3^2y_2=0,\\
M_v^-&:\ \mu_vy_0+\xi_3^2y_1+\xi_3y_2=0.
\end{align*}
At $v=\infty$, where the second Hesse parameter is $0$, put
\[
M_\infty^O=\{y_0=0\},\qquad
M_\infty^+=\{y_1=0\},\qquad
M_\infty^-=\{y_2=0\}.
\]
The corresponding extra divisor classes are
\[
\ell_u^\pm:=\bigl[L_u^\pm\times\mathcal E_{z_*/u}\bigr]
\quad(u\in \Sigma_x),
\qquad
m_v^\pm:=\bigl[\mathcal E_v\times M_v^\pm\bigr]
\quad(v\in \Sigma_y),
\]
with the conventions $z_*/0=\infty$ and $\mathcal E_\infty$ for the fiber at
infinity.
\end{definition}

Note that by \cref{thm:schoen-small-resolutions} and the description of the Neron--Severi groups of the two rational elliptic surfaces and the Picard group of their fiber product, the classes
\[
F,\quad D_x,\quad D_y,\quad
\ell_u^\pm\ (u\in\Sigma_x),\quad
m_v^\pm\ (v\in\Sigma_y)
\]
form a basis of $H^2(W,\mathbb Q)$. 

Also the eight reducible fibers are arranged as
\[
\begin{array}{c|cccccccc}
t&0&1/3&\xi_3/3&\xi_3^2/3&-1/6&-\xi_3/6&-\xi_3^2/6&\infty\\
\hline
\mathcal E_t&I_3&I_3&I_3&I_3&-&-&-&-\\
\mathcal E_{z_*/t}&-&-&-&-&I_3&I_3&I_3&I_3.
\end{array}
\]

\begin{lemma}
\label{lem:gamma-component-action}
The action of $\Gamma$ on the component classes has the following properties.
\begin{enumerate}
\item The automorphism $\phi$ sends $u\mapsto\xi_3u$ and
$v\mapsto\xi_3v$, preserves the labels $O,+,-$, rotates the two triples
\[
\left\{\frac13,\frac{\xi_3}{3},\frac{\xi_3^2}{3}\right\},
\qquad
\left\{-\frac16,-\frac{\xi_3}{6},-\frac{\xi_3^2}{6}\right\},
\]
and fixes individually the chosen classes over $t=0$ and $t=\infty$.
\item The automorphism $\psi$ exchanges the pairs
\[
\left\{0,\frac13\right\},\quad
\left\{\frac{\xi_3}{3},\frac{\xi_3^2}{3}\right\},\quad
\left\{\infty,-\frac16\right\},\quad
\left\{-\frac{\xi_3}{6},-\frac{\xi_3^2}{6}\right\}.
\]
In particular, no chosen extra component class is fixed by $\psi$.
\item The central involution $\psi^2$ stabilizes all eight singular base
values and exchanges $\ell_u^+\leftrightarrow\ell_u^-$ and
$m_v^+\leftrightarrow m_v^-$.  The same exchange occurs on the fixed fibers
$t=0,\infty$ for $\phi\psi^2$ and $\phi^2\psi^2$.
\end{enumerate}
\end{lemma}

\begin{proof}
This follows by substituting \eqref{eq:special-phi-action} and
\eqref{eq:special-psi-action} into the component equations.  For example,
$L_u^O$ is carried by $x_0\mapsto\xi_3x_0$ to
$L_{\xi_3u}^O$.  The remaining equations are checked in the same way.
\end{proof}

\begin{proposition}
\label{prop:gamma-H2-H4-characters}
The representations of $\Gamma$ on $H^2(W,\mathbb C)$ and
$H^4(W,\mathbb C)$ have the same character:
\[
\begin{array}{c|ccccccc}
\text{class} & \mathsf 1 & \mathsf 2 & \mathsf{3A} & \mathsf{3B}
& \mathsf 4 & \mathsf{6A} & \mathsf{6B} \\
\hline
\operatorname{Tr}(g\mid H^2)=\operatorname{Tr}(g\mid H^4)
&19&3&7&7&3&3&3.
\end{array}
\]
Equivalently,
\[
H^2(W,\mathbb C)\cong H^4(W,\mathbb C)
\cong5\rho_1\oplus2\rho_5\oplus2\rho_6\oplus2\rho_7.
\]
\end{proposition}

\begin{proof}
The classes $F,D_x,D_y$ are fixed by every element of $\Gamma$ and
contribute $3$ to every trace.  For $\psi,\psi^2,\phi\psi^2$ and
$\phi^2\psi^2$, no chosen extra component class is fixed.  For $\phi$ and
$\phi^2$, precisely the four classes
$\ell_0^\pm,m_\infty^\pm$ are fixed among the extra components.  Thus the
traces on $H^2$ are $3$ on the classes
$\mathsf 2,\mathsf 4,\mathsf{6A},\mathsf{6B}$ and $7$ on
$\mathsf{3A},\mathsf{3B}$, while the identity trace is $19$.
Averaging an ample class over $\Gamma$ and applying Hard Lefschetz gives a
$\Gamma$-equivariant isomorphism
$H^2(W,\mathbb C)\simeq H^4(W,\mathbb C)$.  Taking inner products with the
character table \eqref{eq:SL2F3-character-table} gives the stated decomposition.
\end{proof}

\begin{theorem}
\label{thm:gamma-H3-decomposition}
The action of $\Gamma$ on the middle cohomology $H^3$ is
\[
H^3(W,\mathbb C)
\cong
2\rho_1\oplus2\rho_2\oplus2\rho_3\oplus4\rho_4
\oplus2\rho_5\oplus2\rho_6\oplus6\rho_7.
\]
\end{theorem}

\begin{proof}
Since $H^1(W)=H^5(W)=0$, the Lefschetz formula gives
\[
\operatorname{Tr}(g\mid H^3)
=
2+\operatorname{Tr}(g\mid H^2)
+\operatorname{Tr}(g\mid H^4)-\Lambda(g).
\]
Using \cref{lem:gamma-fixed-loci,prop:gamma-H2-H4-characters}, the character of
$H^3(W,\mathbb C)$ on the seven conjugacy classes is
\[
(40,\,8,\,-2,\,-2,\,0,\,2,\,2).
\]
Its inner products with the irreducible characters in
\eqref{eq:SL2F3-character-table} are respectively
$2,2,2,4,2,2,6$, which proves the decomposition.
\end{proof}

\begin{proposition}
\label{prop:gamma-invariant-Hodge}
Every element of $\Gamma$ fixes $\Omega_{z_*}$.  Put
\[
V_1:=H^3(W,\mathbb Q)^{\Gamma}.
\]
Consequently $V_1$ has dimension $2$ over $\mathbb Q$ and
\[
V_1\otimes_{\mathbb Q}\mathbb C
=
H^{3,0}(W)\oplus H^{0,3}(W).
\]
\end{proposition}

\begin{proof}
By \cref{thm:gamma-H3-decomposition}, the $\Gamma$-invariant part of
$H^3(W,\mathbb C)$ consists of the two copies of the trivial character and is
therefore two-dimensional.  It remains to identify its Hodge type.
Note that the form $\Omega_{z_*}$ is written as $\omega_t\wedge\omega_{z_*/t}\wedge\frac{dt}{t}$ in \eqref{eq:holomorphic-three-form-Wz}. We will prove that it is fixed by the two generators $\phi,\psi$.

For $\phi$, the numerator and denominator in the residue expression for the
relative one-form on the first factor are both multiplied by $\xi_3$, while
on the second factor they are both multiplied by $\xi_3^2$; moreover
$dt/t$ is fixed.

For $\psi$, let 
\begin{equation}
    \label{eq:Xi}
    \Xi=\begin{pmatrix}1&1&1\\ 1&\xi_3&\xi_3^2\\1&\xi_3^2&\xi_3\end{pmatrix}
\end{equation}
be the first $3\times3$ matrix in~\eqref{eq:special-psi-action}, let $\overline{\Xi}$ be the second one, and write
$F_t=x_0x_1x_2-t(x_0^3+x_1^3+x_2^3)$.  Direct expansion gives
\[
F_{f(t)}(\Xi x)=-\frac9{6t+1}F_t(x),\qquad
\det \Xi=-3\sqrt{-3},
\]
and hence
\[
\psi^*\omega_{f(t)}
=
\frac{\sqrt{-3}}3(6t+1)\omega_t.
\]
For the second factor, one similarly obtains
\[
\psi^*\omega_{f(z_*/t)}
=
-\frac{\sqrt{-3}}3(6z_*/t+1)\omega_{z_*/t}.
\]
Finally,
\[
\frac{df(t)}{f(t)}
=
\frac{9t}{(3t-1)(6t+1)}\frac{dt}{t},
\qquad
6z_*/t+1=\frac{3t-1}{3t}.
\]
The three multipliers cancel, so
$\psi^*\Omega_{z_*}=\Omega_{z_*}$.  
Since $\Omega_z$ spans $H^{3,0}(\mathcal W_z)$ and its conjugate spans $H^{0,3}(\mathcal W_z)$, also pullback commutes with complex conjugation, both $H^{3,0}(W)$ and
$H^{0,3}(W)$ are fixed pointwise by $\Gamma$. 
\end{proof}

\subsection{Splitting of the Hodge structures}
\label{sec:splitting_Hodge}

In this subsection, we study the splitting of the rational Hodge structure. 
Recall that $H^3(W,\mathbb Q)$ carries a pure rational
Hodge structure of weight three, with
\[
H^3(W,\mathbb C)=\bigoplus_{p+q=3}H^{p,q}(W),
\qquad
\overline{H^{p,q}(W)}=H^{q,p}(W).
\]
It is polarized by the alternating form given by
\begin{align}
Q(\alpha,\beta):=-\int_W\alpha\smile\beta.
\label{eq:pol-h3}
\end{align}
The invariant subspace under the action of any finite group
of automorphisms of $W$ is a polarized rational Hodge
substructure, and its orthogonal complement is again
a rational Hodge substructure. The restrictions of $Q$ polarize both Hodge substructures.

\begin{lemma}
\label{lem:G27-invariant-H3}
The rational Hodge structure $H^3(\mathcal W_z,\mathbb Q)^{{G_{27}}}$ has dimension
$4$ and Hodge type
\[
(1,1,1,1).
\]
The Picard Fuchs operator of the three-form $\Omega_z$ is AESZ 4 with $w=z^3$. 
\end{lemma}

\begin{proof}
Take $w=z^3=(3\lambda)^{-6}$. 
By \cref{lem:finite-quotient-cohomology,thm:quotient-H3-comparison},
$H^3(\mathcal W_{z},\mathbb Q)^{G_{27}}\cong
H^3(\mathcal V_{\lambda},\mathbb Q)^{G_{81}}$
as rational Hodge structures. We compute $H^3(\mathcal V_{\lambda},\mathbb Q)^{G_{81}}$
using Griffiths--Dwork reduction. Write
\[
F_x=s_x-3\lambda p_y,\qquad
F_y=s_y-3\lambda p_x
\]
as in \eqref{eq:Vlambda-equations}. Introduce variables $a,b$ and Cayley polynomial $P=aF_x+bF_y$ and the Jacobian ring $R$ of $P$ respects to variables $x_i,y_j,a,b$
\[R={\mathbb C[x_0,x_1,x_2,y_0,y_1,y_2,a,b]}/{(\partial P)}.\]
Consider the variables with bidegree given as follows: 
\[
\deg x_i=\deg y_j=(1,0),
\qquad
\deg a=\deg b=(-3,1).
\]
The Griffiths residue description of a smooth
bicubic complete intersection identifies
\[
H^{3-q,q}(\mathcal V_{\lambda})
\cong R_{(0,q)},\qquad 0\le q\le3.
\]

By the diagonal action \eqref{eq:G81-action} of $G_{81}$, the equations $F_x,F_y$ transform by inverse characters. Let $a,b$ transform by the characters that makes $P$ invariant. 
So the determinant character is trivial, and the residue identification above is $G_{81}$-equivariant. Set $u=ap_y\in R_{(0,1)}$. Note that this is $G_{81}$-invariant. The Jacobian relations are
\[
F_x=F_y=0,\qquad ax_i^2=\lambda bx_jx_k, \qquad by_i^2=\lambda ay_jy_k,
\]
where $\{i,j,k\}=\{0,1,2\}$.
Multiplying the derivative relations by the corresponding coordinates and summing also gives
$u=ap_y=bp_x$ in $R_{(0,1)}$. 
Using these relations and the characters in \eqref{eq:G81-action}, direct monomial reduction gives the following table. 
\[
\begin{array}{c|cccc}
q&0&1&2&3\\ \hline
\dim \mathbb C[x_i,y_j,a,b]_{(0,q)}^{G_{81}}
 &1&8&36&156\\
\dim (\partial P)_{(0,q)}^{G_{81}}
 &0&7&35&155\\ \hline
R_{(0,q)}^{G_{81}}
 &\langle u^0\rangle&\langle u^1\rangle &\langle u^2\rangle&\langle u^3\rangle
\end{array}
\]

Thus we obtain the four $G_{81}$-invariant Hodge components with residue representatives
\[
H^{3,0}: [1],\quad
H^{2,1}: [ap_y],\quad
H^{1,2}: [a^2p_y^2],\quad
H^{0,3}: [a^3p_y^3].
\]

Consequently, $\dim_{\mathbb Q}H^3(\mathcal W_{z},\mathbb Q)^{G_{27}}=4$, and its Hodge type is $(1,1,1,1)$. And the Picard Fuchs statement is either from \cref{prop:holomorphic-period-Wz} or the standard Griffiths--Dwork reduction. 
\end{proof}

\begin{lemma}
\label{lem:G27-preserves-three-form}
Every element of ${G_{27}}$ fixes $\Omega_{z}$.  Consequently,
\[
H^{3,0}(\mathcal W_{z})\oplus H^{0,3}(\mathcal W_{z})\subset
H^3(\mathcal W_{z},\mathbb C)^{{G_{27}}}.
\]
\end{lemma}

\begin{proof}
By \cref{prop:invariant-three-forms}, $h^*\Omega_z=\Omega_z$ for any $h\in G_{27}$. 
Since $\Omega_z$ spans $H^{3,0}(\mathcal W_z)$, its conjugate spans $H^{0,3}(\mathcal W_z)$ and pullback commutes with complex conjugation, both $H^{3,0}(\mathcal W_z)$ and
$H^{0,3}(\mathcal W_z)$ are fixed pointwise by $G_{27}$. 
\end{proof}

\begin{proposition}
\label{prop:rank-two-motivic-splitting}
Define
\[
V_2:=
V_1^\perp\cap H^3(W,\mathbb Q)^{{G_{27}}},
\]
where orthogonality is taken with respect to the polarization \eqref{eq:pol-h3} on $H^3(W)$.
Then
\[
H^3(W,\mathbb Q)^{{G_{27}}}=V_1\oplus V_2.
\]
Both summands have dimension $2$ over $\mathbb Q$, and their Hodge types are
\[
V_1:(1,0,0,1),
\qquad
V_2:(0,1,1,0).
\]
Moreover, $V_1$ is a rank-two submotive of the rank-four
${G_{27}}$-invariant motive; $V_2$ is its orthogonal complement for the
polarization.
\end{proposition}

\begin{proof}
By \cref{prop:gamma-invariant-Hodge,lem:G27-preserves-three-form},
\[
V_1\otimes\mathbb C
=
H^{3,0}(W)\oplus H^{0,3}(W)
\subset
H^3(W,\mathbb C)^{{G_{27}}}.
\]
The invariant Hodge structure on the right is polarized and has Hodge type
$(1,1,1,1)$ by \cref{lem:G27-invariant-H3}.  Its orthogonal complement
to $V_1$ is a two-dimensional rational Hodge structure with one
dimension in $H^{2,1}$ and one in $H^{1,2}$, which proves the assertions of Hodge type. 

For the motivic statement, work first over
$\mathbb Q(\xi_3)$.  Averaging the graphs of the elements of $\Gamma$
and of ${G_{27}}$ gives algebraic idempotent correspondences.  Complex
conjugation sends each of the two groups to itself, so the averages are
Galois invariant and descend to motivated correspondences over $\mathbb Q$.
In fact, if the nontrivial element
$\sigma\in\operatorname{Gal}(\mathbb Q(\xi_3)/\mathbb Q)$ is 
the complex conjugation, then for any $h(\alpha,\beta,\gamma)\in G_{27}$, 
\begin{align}
\label{eq:G27-Galois-conjugation}
\sigma(h(\alpha,\beta,\gamma))=h(-\alpha,-\beta,-\gamma)=h(\alpha,\beta,\gamma)^{-1}\in G_{27}
\end{align}
and
\begin{align}
\label{eq:gamma-Galois-conjugation}
\sigma(\phi)=\phi^{-1},\qquad\sigma(\psi)=\psi^{-1},\qquad \Gamma=\langle\phi,\psi\rangle. 
\end{align}

Their Betti realizations are the two invariant subspaces considered above.
By \cref{lem:G27-preserves-three-form}, the ${G_{27}}$-projector acts as the identity on $V_1$.  
Since the category of pure motivated motives with
rational coefficients is semisimple \cite{andre96}, the
motive $V_1$ is a submotive of rank-four ${G_{27}}$-invariant, and the
polarization gives the complementary motive with realization $V_2$. 
\end{proof}

\begin{proposition}
\label{prop:W-rank-two-attractor}
The special fiber $W=\mathcal W_{-1/18}$ is a rank-two attractor.  More
precisely,
\[
\operatorname{rank}_{\mathbb Z}
\bigl(V_1\cap H^3(W,\mathbb Z)\bigr)=2,
\]
and
\[
\bigl(V_1\cap H^3(W,\mathbb Z)\bigr)\otimes_{\mathbb Z}\mathbb C
=
H^{3,0}(W)\oplus H^{0,3}(W).
\]
\end{proposition}

\begin{proof}
The averaging operator
$|\Gamma|^{-1}\sum_{g\in\Gamma}g^*$ is defined over $\mathbb Q$ on
$H^3(W,\mathbb Q)$, so $V_1$ is a two-dimensional rational subspace of
$H^3(W,\mathbb Q)$.  The image of $H^3(W,\mathbb Z)$ modulo torsion is a
full lattice in $H^3(W,\mathbb Q)$; its intersection with a rational
two-plane is therefore a full lattice of rank $2$ in that plane.  The
second assertion follows from \cref{prop:gamma-invariant-Hodge}.  This is precisely
the rank-two attractor condition.
\end{proof}

\begin{theorem}[Rank-two attractors on the crepant resolutions]
\label{thm:crepant-resolution-attractors}
Suppose that
\[
\lambda_*:=\sqrt{-2},\qquad z_*:=-\frac1{18}=\frac1{9\lambda_*^2},
\qquad w=z_*^3=(3\lambda_*)^{-6}=-\frac1{5832}.
\]
Let
\[
\rho_X:\widetilde{\mathcal X}_{\lambda_*}\longrightarrow\mathcal X_{\lambda_*},\qquad\rho_Y:\widetilde{\mathcal Y}_{z_*}\longrightarrow\mathcal Y_{z_*}
\]
be projective crepant resolutions. Then the two resolutions are smooth
projective Calabi--Yau threefolds, and both pullbacks
\[
\begin{aligned}
\rho_X^*:\ H^3(\mathcal X_{\lambda_*},\mathbb Q)
&\xrightarrow{\sim} H^3(\widetilde{\mathcal X}_{\lambda_*},\mathbb Q),\\
\rho_Y^*:\ H^3(\mathcal Y_{z_*},\mathbb Q)
&\xrightarrow{\sim} H^3(\widetilde{\mathcal Y}_{z_*},\mathbb Q)
\end{aligned}
\]
are isomorphisms of rational Hodge structures. In particular, both
$\widetilde{\mathcal X}_{\lambda_*}$ and $\widetilde{\mathcal Y}_{z_*}$
are rank-two attractors.
\end{theorem}

\begin{proof}
Put $V:=\mathcal V_{\lambda_*}, X:=\mathcal X_{\lambda_*}$ and $W:=\mathcal W_{z_*}, Y:=\mathcal Y_{z_*}$ for simplicity, and write
$q_X:V\to X$ and
$q_Y:W\to Y$ for the quotient maps.
By \cref{lem:finite-quotient-cohomology}, pullback gives isomorphisms
of pure rational Hodge structures
\[
\begin{aligned}
H^3(X,\mathbb Q)
&\xrightarrow[\sim]{q_X^*}H^3(V,\mathbb Q)^{G_{81}},\\
H^3(Y,\mathbb Q)
&\xrightarrow[\sim]{q_Y^*}H^3(W,\mathbb Q)^{G_{27}}.
\end{aligned}
\]
Moreover, \cref{thm:quotient-H3-comparison} identifies the two quotient
cohomology groups through $\Phi_{\lambda_*}$. Hence
\cref{lem:G27-invariant-H3} implies both of them have dimension four: 
\begin{equation}\label{eq:crepant-attractor-quotient-b3}
\dim_{\mathbb Q}H^3(X,\mathbb Q)=\dim_{\mathbb Q}H^3(Y,\mathbb Q)=4.
\end{equation}

We first claim that pullback on $H^3$ is injective for either resolution.
More generally, let $T$ ($T=V$ or $W$ in our situation) be a smooth projective threefold with an
algebraic action of a finite group $F$ ($F=G_{81}$ or $G_{27}$ respectively), let $S=T/F$, and let
$\rho:\widetilde S\to S$ be a projective resolution. The cup-product
pairing on $H^3(T,\mathbb Q)$ is nondegenerate and $F$-invariant,
so its restriction to $H^3(T,\mathbb Q)^F$ is nondegenerate.
By \cref{lem:finite-quotient-cohomology}, the formula
\[
(\alpha,\beta)_S
:=\frac1{|F|}\int_T q^*\alpha\smile q^*\beta,
\qquad \alpha,\beta\in H^3(S,\mathbb Q),
\]
defines a nondegenerate pairing on $H^3(S,\mathbb Q)$, where
$q:T\to S$ is the quotient map. Since $q$ has degree $|F|$
and $\rho$ has degree one, the projection formula gives
\[
\int_{\widetilde S}\rho^*\alpha\smile\rho^*\beta
=\frac1{|F|}\int_T q^*\alpha\smile q^*\beta
=(\alpha,\beta)_S.
\]
Consequently $\rho^*:H^3(S,\mathbb Q)\to
H^3(\widetilde S,\mathbb Q)$ is injective.

For the standard Libgober--Teitelbaum crepant resolution of
$X$ one has $h^{2,1}=1$ \cite{libgober1993}.
The two resolutions in the statement and this standard resolution are
birational smooth projective Calabi--Yau threefolds, by
\cref{prop:birational-quotients}. Birational invariance of Betti
numbers for smooth projective Calabi--Yau threefolds~\cite{Batyrev:1999ab} (also obtainable
from their factorization into flops \cite{kawamata08}) therefore gives
\[
\dim_{\mathbb Q}H^3(\widetilde{X},\mathbb Q)
=\dim_{\mathbb Q}H^3(\widetilde{Y},\mathbb Q)
=2+2h^{2,1}=4.
\]
Together with \eqref{eq:crepant-attractor-quotient-b3} and the injectivity
just proved, this shows that $\rho_X^*$ and
$\rho_Y^*$ are isomorphisms. Since pullbacks are morphisms
of mixed Hodge structures and all four cohomology groups are pure of
weight three, they are isomorphisms of rational Hodge structures.

We finally transfer the rank-two decomposition. Set
\[
\begin{aligned}
\Theta_X
&:=q_Y^*\circ\Phi_{\lambda_*}
\circ(\rho_X^*)^{-1}:
H^3(\widetilde{X},\mathbb Q)
\xrightarrow{\sim}H^3(W,\mathbb Q)^{G_{27}}, \\
\Theta_Y
&:=q_Y^*\circ(\rho_Y^*)^{-1}:
H^3(\widetilde{Y},\mathbb Q)
\xrightarrow{\sim}H^3(W,\mathbb Q)^{G_{27}}. 
\end{aligned}
\]
For $S=X$ or $Y$, denote the corresponding
isomorphism by $\Theta_{S}$, and put
\[
V_{i,\widetilde S}:=\Theta_{S}^{-1}(V_i)
\qquad(i=1,2).
\]
By \cref{prop:rank-two-motivic-splitting},
\[
H^3(\widetilde S,\mathbb Q)
=V_{1,\widetilde S}\oplus V_{2,\widetilde S},
\qquad
V_{1,\widetilde S}\otimes_{\mathbb Q}\mathbb C
=H^{3,0}(\widetilde S)\oplus H^{0,3}(\widetilde S).
\]
The two summands are rational Hodge structures of dimension two.
The intersection of the rational two-plane $V_{1,\widetilde S}$
with $H^3(\widetilde S,\mathbb Z)$ modulo torsion is a lattice of
rank two. Thus
\[
\operatorname{rank}_{\mathbb Z}
\bigl(V_{1,\widetilde S}\cap H^3(\widetilde S,\mathbb Z)\bigr)=2,
\]
which proves the rank-two attractor assertion for both resolutions.
\end{proof}

\begin{remark}
Note that for any $u,v$ such that $u^6=v^3=1$, 
we have actually isomorphisms $\mathcal X_{u\lambda_*}\cong\mathcal X_{\lambda_*}$ and $\mathcal Y_{vz_*}\cong\mathcal Y_{z_*}$. 
So the resolved Calabi--Yau threefolds $\widetilde{X}_{u\lambda_*}$ and $\widetilde{Y}_{vz_*}$ are rank-two attractors too. 
\end{remark}

\begin{remark}
It is natural to ask whether we could define our extra automorphism $\psi$ on $\cV_{\lambda_*},\cX_{\lambda_*},\cY_{z_*}$ as well. 
Actually the automorphism $\psi_{\cW}:=\psi$ of $\cW_{z_*}$ in
\eqref{eq:special-psi-action}, has a linear lift to the bicubic $\cV_{\lambda_*}$:
\[
\psi_{\cV}:\cV_{\lambda_*}\longrightarrow\cV_{\lambda_*},
\qquad [x:y]\longmapsto[\Xi x:-\overline{\Xi}y],
\]
where $\Xi$ and $\overline{\Xi}$ were introduced in~\eqref{eq:Xi}. Indeed, with $F_x=\sx-3\lambda_*\py$ and
$F_y=\sy-3\lambda_*\px$, one has
\[
\begin{pmatrix}F_x\circ\psi_{\cV}\\ F_y\circ\psi_{\cV}\end{pmatrix}
=\begin{pmatrix}3&3\lambda_*\\-3\lambda_*&-3\end{pmatrix}
\begin{pmatrix}F_x\\F_y\end{pmatrix},
\qquad \lambda_*^2=-2.
\]
Moreover, $\psi_{\cV}^2$ exchanges $x_1,x_2$ and $y_1,y_2$,
so that $\psi_{\cV}$ is a one-to-one automorphism of order four.
The degree-three rational map \eqref{eq:V-to-W-rational-map}
intertwines $\psi_{\cV}$ and $\psi_{\cW}$.
Its three generic preimages form one orbit of the kernel
$K\cong\mZ_3$ in \eqref{eq:K}; since
$\psi_{\cV}$ commutes with $K$, it acts bijectively on this orbit. 

However, neither automorphism descends to a single-valued map on the respective
quotient. Conjugation by the matrices $\Xi,\overline{\Xi}$ takes nontrivial diagonal
transformations to non-diagonal transformations; on the base, conjugating
$t\mapsto\xi_3^\gamma t$ by $f$ does not preserve the scaling subgroup
when $\gamma\ne0$. Thus
\[
G_{27}\cap\psi_{\cW}G_{27}\psi_{\cW}^{-1}=\{1\},\qquad
G_{81}\cap\psi_{\cV}G_{81}\psi_{\cV}^{-1}=K.
\]
Consequently their graph images in
$\cY_{z_*}\times\cY_{z_*}$ and
$\cX_{\lambda_*}\times\cX_{\lambda_*}$ are irreducible algebraic
correspondences of generic bidegree $(27,27)$, since
$|G_{27}|=27$ and $|G_{81}|/|K|=27$. So we can think of them as 27-to-27 correspondences. 

These two correspondences are explicit on the open dense subset
$\cT_{z_*^3}^{\circ}$ of \cref{prop:Y2T}.
Given coordinates $(u_i,v_i)$, choose $t,x_i,y_i$ satisfying
\[
t^3=u_0u_1u_2,\qquad x_i^3=u_i,\quad \px=t,\qquad
y_i^3=v_i,\quad \py=z_*/t.
\]
The corresponding image points $(u_i',v_i')$ are
\[
u_i'=\frac{\bigl(\sum_{j=0}^2\xi_3^{ij}x_j\bigr)^3}
{3(1+6t)},\qquad
v_i'=\frac{\bigl(\sum_{j=0}^2\xi_3^{-ij}y_j\bigr)^3}
{3(1+6z_*/t)},\qquad 0\le i\le2.
\]
There are three choices of $t$ and, for each, three projective choices
of each of $x$ and $y$, giving $3\cdot3\cdot3=27$ distinct generic
images. Although $\psi_{\cV}^4=\psi_{\cW}^4=1$, the fourth iterate
of either quotient correspondence is not the diagonal trivial correspondence. 

Finally, by \cref{lem:finite-quotient-cohomology,thm:quotient-H3-comparison},
both quotient cohomology groups are identified with
$H^3(\cW_{z_*},\mQ)^{G_{27}}=V_1\oplus V_2$.
The normalized action of the correspondence is
$e_{G_{27}}\psi_{\cW}^*$ on this invariant subspace, where
$e_{G_{27}}=27^{-1}\sum_{h\in G_{27}}h^*$.
A fixed-locus computation gives
\[
\frac1{27}\sum_{h\in G_{27}}\Lambda(\psi_{\cW}h)
=\frac{72+54+54}{27}=\frac{20}{3},
\qquad
\operatorname{Tr}\bigl(e_{G_{27}}\psi_{\cW}^*\mid H^3\bigr)
=1+3+3+1-\frac{20}{3}=\frac43.
\]
Here $72,54,54$ are the sums over $\alpha,\beta\in\mZ_3$
for $\gamma=0,1,2$ in \eqref{eq:G27-action}, respectively.
As $\psi_{\cW}$ fixes $V_1$ pointwise, and $\psi_{\cW}^2$
acts trivially on the four-dimensional invariant cohomology
(it fixes the holomorphic three-form and its Gauss--Manin derivatives),
the normalized action is self-adjoint for the polarization.
Its restriction to the two-dimensional symplectic space $V_2$ is
therefore scalar; the trace determines
\[
\left.e_{G_{27}}\psi_{\cW}^*\right|_{V_1}=1,
\qquad
\left.e_{G_{27}}\psi_{\cW}^*\right|_{V_2}=-\frac13.
\]
These eigenvalues apply to $H^3(\cX_{\lambda_*},\mQ)$ as well,
through $\Phi_{\lambda_*}$. 
Thus the correspondence does not induce an automorphism preserving
the integral lattice or the polarization (and its fourth
power is not the identity). 
\end{remark}

\section{Modularity for the AESZ 4 attractor}
\label{sec:point-counting-modularity}

\subsection{Arithmetic realization}
\label{sec:Arith_real}

We pass from the polarized rational Hodge structures in the preceding
section to their arithmetic realizations to complete the modularity argument.  
This transition uses algebraic correspondences and the Betti--\'{e}tale comparison isomorphism as main tools.  

We first state our Frobenius and Tate-twist conventions. 
Throughout this section, $G_{\mathbb Q}=\operatorname{Gal}(\overline{\mathbb Q}/\mathbb Q)$, and $\chi_\ell:G_{\mathbb Q}\to\mathbb Z_\ell^\times$ denotes the $\ell$-adic cyclotomic character describing the Galois action on $\mathbb Q_\ell(1)$. We use the convention $\mathbb Q_\ell(n)=\mathbb Q_\ell(1)^{\otimes n}$, and $\mathbb Q_\ell(n)$ has Galois character $\chi_\ell^n$. 

For a prime $p\ne\ell$ at which the representations under consideration are unramified, $\operatorname{Frob}_p$ denotes \emph{geometric Frobenius}, the inverse of arithmetic Frobenius. Thus
\[
\chi_\ell(\operatorname{Frob}_p)=p^{-1},
\qquad
\operatorname{Frob}_p\bigm|_{\mathbb Q_\ell(-n)}=p^n.
\]

For a normalized newform $f=\sum_{n\ge1}a_n(f)q^n$ of weight $k$ and trivial nebentypus, we use the its associated $\ell$-adic Galois representation $\rho_{f,\ell}$, characterized at primes $p\nmid N\ell$ by
\[
\det(1-T\operatorname{Frob}_p\mid\rho_{f,\ell})
=1-a_p(f)T+p^{k-1}T^2.
\]
In particular, under these conventions,
\[
\det\rho_{f,\ell}=\chi_\ell^{-(k-1)}.
\]

Next, we review and apply the Betti--\'{e}tale comparison theorem in our context.
Although $W$ is defined over $\mathbb Q$, the automorphism elements
of $G_{27}$ and $\Gamma$ are generally defined only over $\mathbb Q(\xi_3)$.
However, the Galois-stability of each group (check \eqref{eq:G27-Galois-conjugation} and \eqref{eq:gamma-Galois-conjugation} in the proof of \cref{prop:rank-two-motivic-splitting}) 
makes its average into a projector defined over $\mathbb Q$.  More explicitly,
for every prime $\ell$ set
\begin{align*}
 e_{G_{27}}&=\frac1{|G_{27}|}\sum_{h\in G_{27}}h^*,
 &e_{\Gamma}&=\frac1{|\Gamma|}\sum_{g\in \Gamma}g^*,\\
 U_{G_{27},\ell}&=e_{G_{27}}H^3_{\mathrm{\acute et}}
       \left(W_{\overline{\mathbb Q}},\mathbb Q_\ell\right)
 &U_{\Gamma,\ell}&=e_{\Gamma}H^3_{\mathrm{\acute et}}
       \left(W_{\overline{\mathbb Q}},\mathbb Q_\ell\right). 
\end{align*}
For $\sigma\in G_{\mathbb Q}$ one has
\[
 \sigma e_{G_{27}}\sigma^{-1}
   =\frac1{|G_{27}|}\sum_{h\in G_{27}}(h^\sigma)^*=e_{G_{27}},\qquad
   \sigma e_{\Gamma}\sigma^{-1}
   =\frac1{|\Gamma|}\sum_{g\in\Gamma}(g^\sigma)^*=e_{\Gamma},
\]
so $U_{G_{27},\ell},U_{\Gamma,\ell}$ are $G_{\mathbb Q}$-subrepresentations.  
Here all graph averages have rational coefficients.

We now explain how the cohomological decomposition
obtained in the preceding section passes to $\ell$-adic cohomology. 
Fix an embedding $\iota:\overline{\mathbb Q}\hookrightarrow\mathbb C$.
For every prime $\ell$, the Betti--\'etale comparison theorem
gives a canonical isomorphism relative to this embedding,
\[
\operatorname{comp}_{\ell}:
H^3_{\mathrm B}(W(\mathbb C),\mathbb Q)
\otimes_{\mathbb Q}\mathbb Q_\ell
\xrightarrow{\ \sim\ }
H^3_{\mathrm{\acute et}}
(W_{\overline{\mathbb Q}},\mathbb Q_\ell).
\]
This comparison is functorial for algebraic morphisms and
compatible with algebraic cycle classes, cup products, and
Poincar\'e duality; for details see
\cite[Secs.~21, 23, and 24]{milneLEC}.

Let
\[
e_{G_{27}}^\mathrm B=\frac1{|G_{27}|}\sum_{h\in G_{27}}h^*,\qquad
 e_{\Gamma}^\mathrm B=\frac1{|\Gamma|}\sum_{g\in \Gamma}g^*,
\]
denote the averaging operator on Betti cohomology.
Functoriality of the comparison isomorphism therefore gives
\[
\operatorname{comp}_{\ell}\circ
\left(e_{G_{27}}^{\mathrm B}\otimes 1\right)
=e_{G_{27}}\circ\operatorname{comp}_{\ell},\qquad
\operatorname{comp}_{\ell}\circ
\left(e_\Gamma^{\mathrm B}\otimes 1\right)
=e_\Gamma\circ\operatorname{comp}_{\ell}.
\]
Consequently,
\[
U_{G_{27},\ell}
\cong
H^3_{\mathrm B}(W(\mathbb C),\mathbb Q)^{G_{27}}
\otimes_{\mathbb Q}\mathbb Q_\ell,\qquad
U_{\Gamma,\ell}
\cong
H^3_{\mathrm B}(W(\mathbb C),\mathbb Q)^{\Gamma}
\otimes_{\mathbb Q}\mathbb Q_\ell
\]
as $\mathbb Q_\ell$-vector spaces, compatibly with the
averaging projectors.

By the dimension results of the preceding section,
\[\begin{aligned}
\dim_{\mathbb Q}
H^3_{\mathrm B}(W(\mathbb C),\mathbb Q)^\Gamma=2,&\quad
\dim_{\mathbb Q}
H^3_{\mathrm B}(W(\mathbb C),\mathbb Q)^{G_{27}}=4,
\\
H^3_{\mathrm B}(W(\mathbb C),\mathbb Q)^\Gamma
&\subseteq
H^3_{\mathrm B}(W(\mathbb C),\mathbb Q)^{G_{27}}.
\end{aligned}\]
Since the comparison isomorphism intertwines the
projectors, these statements imply
\[
\dim_{\mathbb Q_\ell}U_{\Gamma,\ell}=2,
\qquad
\dim_{\mathbb Q_\ell}U_{G_{27},\ell}=4,
\qquad
U_{\Gamma,\ell}\subseteq U_{G_{27},\ell}.
\]

We next verify that the orthogonal decomposition in~\ref{prop:rank-two-motivic-splitting} is
preserved. Firstly the $\ell$-adic Poincar\'e pairing
\[
\langle\ ,\ \rangle_\ell:
H^3_{\mathrm{\acute et}}
(W_{\overline{\mathbb Q}},\mathbb Q_\ell)
\times
H^3_{\mathrm{\acute et}}
(W_{\overline{\mathbb Q}},\mathbb Q_\ell)
\longrightarrow
\mathbb Q_\ell(-3)
\]
is nondegenerate and $G_{\mathbb Q}$-equivariant.
It corresponds to the rational cup-product pairing used to define
the polarization in the preceding section under Betti--\'etale comparison. Note that every geometric automorphism $g$ preserves this pairing, and hence
\[
\langle g^*x,y\rangle_\ell=\langle x,(g^{-1})^*y\rangle_\ell.
\]
Averaging over the finite group shows that
\[
\langle e_{G_{27}}x,y\rangle_\ell=\langle x,e_{G_{27}} y\rangle_\ell, \qquad
\langle e_{\Gamma}x,y\rangle_\ell=\langle x,e_\Gamma y\rangle_\ell. 
\]
Thus $e_{G_{27}},e_{\Gamma}$ are self-adjoint idempotent
for the Poincar\'e pairing, and their images are
nondegenerate. In particular, the inclusion
$U_{\Gamma,\ell}\subseteq U_{G_{27},\ell}$
implies
\[
e_{G_{27}}e_\Gamma=e_\Gamma e_{G_{27}}=e_\Gamma
\]
on third $\ell$-adic cohomology.
It follows that $e_{G_{27}}-e_\Gamma$
is a self-adjoint idempotent, whose image is precisely
the orthogonal complement of $U_{\Gamma,\ell}$
inside $U_{G_{27},\ell}$.

We may therefore define
\[
V_{1,\ell}:=U_{\Gamma,\ell}, \qquad
V_{2,\ell}:=\operatorname{im}(e_{G_{27}}-e_\Gamma)
  =U_{G_{27},\ell}\cap V_{1,\ell}^{\perp}.
\]
Then
\[
U_{G_{27},\ell}
=
V_{1,\ell}\oplus V_{2,\ell},
\qquad
\dim_{\mathbb Q_\ell}V_{1,\ell}
=
\dim_{\mathbb Q_\ell}V_{2,\ell}=2.
\]

Finally, this decomposition is defined over $\mathbb Q$
in the arithmetic sense. Indeed, as established above,
both averaging projectors commute with
$G_{\mathbb Q}$ because the corresponding groups
are stable under Galois conjugation.
Therefore $V_{1,\ell}$ and $V_{2,\ell}$
are $G_{\mathbb Q}$-stable subrepresentations. 
These two rank-two $\ell$-adic realizations are thus
the arithmetic realizations of the corresponding
motivated direct summands constructed in
\cref{prop:rank-two-motivic-splitting}; see also
\cite{andre96,scholl94}.

\begin{lemma}
\label{lem:good-reduction-W}
The threefold $W$ has a smooth projective model over
$\mathbb Z[1/6]$.  For \emph{every} prime $\ell$ (including $\ell=2$ we will use later)
and every prime $p>3$ distinct from $\ell$, the representations on
$V_{1,\ell}$ and $V_{2,\ell}$ are unramified at $p$.
\end{lemma}

\begin{proof}
Write $t=t_0/t_1$ in homogeneous coordinates
$(t_0:t_1)$ on the base $\mathbb P^1$.
The fiber product is given in
$\mathbb P^2\times\mathbb P^2\times\mathbb P^1$ by
\[
t_0(x_0^3+x_1^3+x_2^3)-t_1x_0x_1x_2=0,
\qquad
t_1(y_0^3+y_1^3+y_2^3)+18t_0y_0y_1y_2=0.
\]
The two Hesse parameters are $t$ and
$z_*/t=-t_1/(18t_0)$, respectively.

For $p>3$, base change to $\overline{\mathbb F}_p$, 
the singular fibers of the first surface occur at
\[
t=0,\quad \frac13,\quad \frac{\xi_3}{3},
\quad \frac{\xi_3^2}{3},
\]
and those of the second at
\[
t=\infty,\quad -\frac16,\quad
-\frac{\xi_3}{6},\quad -\frac{\xi_3^2}{6}.
\]
Since equality of two finite values would give
\[\frac1{27}\ne -\frac1{216},\qquad \infty\ne\frac1{27},\qquad0\ne-\frac1{216},\qquad0\ne\infty\]
after cubing, that is bt $9\ne 0$ and $2\ne 0$ in $\mathbb F_p$, hence these sets are disjoint. 
As both rational elliptic surfaces are smooth,
their fiber product is smooth over $\mathbb Z[1/6]$.
Also it is projective by construction.

The actions of $G_{27}$ and $\Gamma$ extend over $\mathbb Z[\xi_3,1/6]$. Indeed, the diagonal factors are units, the Fourier matrices are invertible after inverting $3$, and the base transformation $t\mapsto(1-3t)/(3+18t)$ has determinant $-27$. 

Thus the group actions reduce modulo every $p>3$, and the corresponding twisted Frobenius fixed-point counts are well-defined. The unramifiedness assertion follows from smooth proper base change.
\end{proof}

\begin{lemma}
\label{lem:G27-even-cohomology}
For every prime $\ell$,
\[
H^2_{\mathrm{\acute et}}(W_{\overline{\mathbb Q}},\mathbb Q_\ell)^{{G_{27}}}
\cong\mathbb Q_\ell(-1)^{\oplus7},
\qquad
H^4_{\mathrm{\acute et}}(W_{\overline{\mathbb Q}},\mathbb Q_\ell)^{{G_{27}}}
\cong\mathbb Q_\ell(-2)^{\oplus7}.
\]
Similarly,
\[
H^2_{\mathrm{\acute et}}(W_{\overline{\mathbb Q}},\mathbb Q_\ell)^{\Gamma}
\cong\mathbb Q_\ell(-1)^{\oplus5},
\qquad
H^4_{\mathrm{\acute et}}(W_{\overline{\mathbb Q}},\mathbb Q_\ell)^{\Gamma}
\cong\mathbb Q_\ell(-2)^{\oplus5}.
\]
\end{lemma}

\begin{proof}

We compute the invariant parts of $H^2$ by giving explicit
bases of divisor classes defined over $\mathbb Q$.
Their $\ell$-adic cycle classes will then give copies
of $\mathbb Q_\ell(-1)$.

Recall that each Hesse surface is the blow-up of
$\mathbb P^2$ at the nine base points of the Hesse
pencil. Its Picard group is generated by the
hyperplane class and the nine exceptional sections.
By \cite{sch88} and Section~\ref{sec:hesse-fiber-products},
the divisor group of $W$ is obtained from the two
rank-ten Picard groups by identifying their common
fiber class $F$.

Let $H_x$ and $H_y$ denote the pullbacks of the two
hyperplane classes. On each Hesse surface, $G_{27}$
has three orbits of size three on the nine base points,
one on each coordinate line. The invariant Picard
rank is therefore $1+3=4$ on each surface, and
$4+4-1=7$ on $W$. A basis is explicitly given by
\[
\bigl\{F,\quad H_x, H_y,\quad
\ell_0^+,\ell_0^-,\quad
m_\infty^+,m_\infty^-\bigr\}.
\]
These are all $G_{27}$-invariant divisor classes
defined over $\mathbb Q$.

For $\Gamma$, the zero section on each surface is fixed,
while the other eight flex points form a single orbit
under the natural $\operatorname{SL}(2,\mathbb F_3)$ action.
Thus the invariant rank is $1+2=3$ on each surface,
and $3+3-1=5$ on $W$, with basis
\[
\bigl\{F,\quad H_x, H_y,\quad D_x,D_y\bigr\}.
\]
Again, all five classes are defined over $\mathbb Q$.
The independence of both bases follows from the
rank-$19$ divisor presentation.

Since $h^{2,0}(W)=0$, all of $H^2$ is algebraic.
By the cycle class map and Betti--\'etale comparison,
each of the above rational divisor classes spans a
copy of $\mathbb Q_\ell(-1)$.
Hence
\[
H^2_{\mathrm{\acute et}}(W_{\overline{\mathbb Q}},
\mathbb Q_\ell)^{G_{27}}
\simeq \mathbb Q_\ell(-1)^{\oplus 7},
\qquad
H^2_{\mathrm{\acute et}}(W_{\overline{\mathbb Q}},
\mathbb Q_\ell)^\Gamma
\simeq \mathbb Q_\ell(-1)^{\oplus 5}
\]
as $G_{\mathbb Q}$-representations for every $\ell$.

Finally, Hard Lefschetz gives
$H^2\simeq H^4(1)$ by cup product with
$c_1(\mathcal O(1,1,1)|_W)$.
This ample class is defined over $\mathbb Q$ and
invariant under both groups, so the isomorphism is
equivariant for $G_{\mathbb Q}$, $G_{27}$, and
$\Gamma$. The assertions for $H^4$ follow.
\end{proof}

\subsection{Point counting and trace comparisons}
\label{sec:point_count}

We first compute the Frobenius traces at good primes of the two rank-two Galois representations by equivariant point counting. 
Then we identify these traces with coefficients of appropriate newforms in the modular forms database. 
Finally we apply the Serre--Faltings--Livn\'{e}--Sch\"{u}tt criterion to prove modularity.

For a good prime power $q=p^r$ with $p>3$, choose an algebraic
closure of $\mathbb F_p$ containing a primitive cube root of unity,
and reduce the automorphisms of both groups to the geometric special
fiber.  They are defined over $\mathbb F_{p^2}$ when $p\equiv2\pmod3$
(and over $\mathbb F_p$ when $p\equiv1\pmod3$), but need not be
defined over $\mathbb F_p$.  Accordingly, the \emph{twisted Frobenius}
$g\circ\operatorname{Frob}_q$ is considered as an endomorphism over
$\overline{\mathbb F}_p$; it need not commute with Frobenius.
Grothendieck--Lefschetz applies to this endomorphism, and its fixed
points are simple because Frobenius has zero differential.  The
projector average commutes with Frobenius, because
Frobenius conjugates each group into itself. 

For an automorphism element $g$ in either ${G_{27}}$ or $\Gamma$, put
\[
\Lambda_q(g):=
\sum_{i=0}^6(-1)^i
\operatorname{Tr}\!\left(
g^*\circ\operatorname{Frob}_q^*
\mid
H^i_{\mathrm{\acute et}}
(W_{\overline{\mathbb F}_p},\mathbb Q_\ell)
\right).
\]
The Lefschetz fixed-point formula identifies $\Lambda_q(g)$ with the number
of geometric fixed points of $g\circ\operatorname{Frob}_q$.  Since
\begin{align*}
\operatorname{Tr}(\operatorname{Frob}_p\,e_\Gamma|H^i)
&=|\Gamma|^{-1}\sum_{g\in \Gamma}\operatorname{Tr}(\operatorname{Frob}_p g^*|H^i), \\
\operatorname{Tr}(\operatorname{Frob}_p\,e_{G_{27}}|H^i)
&=|G_{27}|^{-1}\sum_{h\in G_{27}}\operatorname{Tr}(\operatorname{Frob}_p h^*|H^i)
\end{align*}
and $H^1=H^5=0$, averaging and using
\cref{lem:G27-even-cohomology} gives, for every good prime $p$,
\begin{align}
\label{eq:G27-invariant-Frobenius-trace}
\operatorname{Tr}\bigl(\operatorname{Frob}_p\mid
H^3(W)^{{G_{27}}}\bigr)
&=
1+7p+7p^2+p^3
-\frac1{27}\sum_{h\in{G_{27}}}\Lambda_p(h),\\
\label{eq:V1-Frobenius-trace}
\operatorname{Tr}(\operatorname{Frob}_p\mid V_1)
&=
1+5p+5p^2+p^3
-\frac1{24}\sum_{g\in\Gamma}\Lambda_p(g).
\end{align}
Since $H^3(W)^{{G_{27}}}=V_1\oplus V_2$, the trace on $V_2$ is the
difference of these two traces.

The point-counting computation program we used is a PARI/GP~\cite{PARI2} script in appendix~\ref{app:pari}.  
For the finite primes set
\[
T=\{5,7,11,13,17,19,23,31,37\},
\]
the ${G_{27}}$-invariant and the $\Gamma$-invariant computations give
\[
\begin{array}{c|rrrrrrrrr}
p&5&7&11&13&17&19&23&31&37\\ \hline
\frac1{27}\sum_{h\in{G_{27}}}\Lambda_p(h)
&318&714&2346&3504&7128&9588&15996&36534&60420\\
\operatorname{Tr}(\operatorname{Frob}_p\mid H^3(W)^{{G_{27}}})
&18&22&-90&-32&-72&-68&36&202&76\\
\operatorname{Tr}(\operatorname{Frob}_p\mid V_1)
&3&29&-57&20&-72&-106&174&47&2\\
\operatorname{Tr}(\operatorname{Frob}_p\mid V_2)
&15&-7&-33&-52&0&38&-138&155&74\\
\operatorname{Tr}(\operatorname{Frob}_p\mid V_2)/p
&3&-1&-3&-4&0&2&-6&5&2.
\end{array}
\]

\begin{theorem}[Serre--Faltings--Livn\'e--Sch\"utt criterion for bad primes $S\subset\{2,3\}$]
\label{thm:S23-comparison-criterion}
Let $\rho_1$ and $\rho_2$ be two-dimensional $2$-adic Galois
representations of $G_{\mathbb Q}$, unramified outside $\{2,3\}$, with the
same determinant.  Assume that the trace is even at
$\operatorname{Frob}_{11}$ or at $\operatorname{Frob}_{13}$.  Then equality
of the traces of $\rho_1$ and $\rho_2$ at every prime in
\[
T=\{5,7,11,13,17,19,23,31,37\}
\]
implies
$\rho_1^{\mathrm{ss}}\cong\rho_2^{\mathrm{ss}}$.
\end{theorem}

\begin{proof}
This is the finite test in \cite[Proposition~1]{schutt06}, which combines
Serre's uneven-trace argument with Livn\'e's criterion
\cite[Theorem~4.3]{liv87}.
\end{proof}

\begin{theorem}[Modularity of $V_1$ and $V_2$]
\label{thm:V1-V2-modularity}
Let
\[
f=f_{\mathsf{54.4.a.c}}\in S_4(\Gamma_0(54)),
\qquad
g=f_{\mathsf{54.2.a.a}}\in S_2(\Gamma_0(54))
\]
be the normalized newforms 
\cite[\href{https://www.lmfdb.org/ModularForm/GL2/Q/holomorphic/54/4/a/c/}{$f_{\mathsf{54.4.a.c}}$}]{lmfdb}, \cite[\href{https://www.lmfdb.org/ModularForm/GL2/Q/holomorphic/54/2/a/a/}{$f_{\mathsf{54.2.a.a}}$}]{lmfdb} with rational coefficients and trivial character.  Then
\[
\rho_{V_1,2}^{\mathrm{ss}}
\cong
\rho_{f,2}^{\mathrm{ss}},
\qquad
\rho_{V_2,2}^{\mathrm{ss}}
\cong
\rho_{g,2}(-1)^{\mathrm{ss}}.
\]
\end{theorem}

\begin{proof}
By \cref{lem:good-reduction-W}, the two
$\mathbb Q_2$-representations $V_{1,2}$ and $V_{2,2}$ are unramified
at every $p>3$.  Their two-dimensional $G_{\mathbb Q}$-stable
subspaces are obtained using rational projectors. 
A stable $\mathbb Z_2$-lattice exists by intersecting each subspace
with a stable lattice of the ambient proper smooth cohomology.
The modular-form representations are also unramified outside the
primes dividing $2\cdot54$ \cite{deligne71}.

The cup-product pairing on $H^3$ is a $G_{\mathbb Q}$-equivariant
alternating pairing with values in $\mathbb Q_2(-3)$.
Since finite-group averaging is orthogonal for this pairing and
$V_{1,2}\subset H^3(W)^{{G_{27}}}$, its restriction to each of the
two rank-two summands is nondegenerate.  In dimension two this
identifies $\bigwedge^2V_{i,2}$ with $\mathbb Q_2(-3)$,
so the whole determinant characters satisfy
\[
 \det V_{1,2}=\det V_{2,2}=\chi_2^{-3}. 
\]
The newforms have trivial nebentypus, whence
$\det\rho_{f,2}=\chi_2^{-3}$ and
$\det\rho_{g,2}(-1)=\chi_2^{-3}$.
Thus the determinant condition of \cref{thm:S23-comparison-criterion} holds.

The newform coefficients are
\[
\begin{array}{c|rrrrrrrrr}
p&5&7&11&13&17&19&23&31&37\\ \hline
a_p(f)&3&29&-57&20&-72&-106&174&47&2\\
a_p(g)&3&-1&-3&-4&0&2&-6&5&2.
\end{array}
\]
Thus the point-counted traces satisfy
$\operatorname{Tr}(\operatorname{Frob}_p\mid V_1)=a_p(f)$ and
$\operatorname{Tr}(\operatorname{Frob}_p\mid V_2)=p\,a_p(g)$ for every
$p\in T$.  At $p=13$ the two traces $20$ and $-52$ are even.  Applying
\cref{thm:S23-comparison-criterion} to the two pairs gives the asserted
semisimplified isomorphisms.  Equality of the good Euler factors follows
from equality of traces and determinants.
\end{proof}

Now we realize everything in a motivic way to finish the comparison using $\mathbb Q$-models: 
\[
 X^{(0)}:=V^{(0)}/G_{81},\qquad
 V^{(0)}:=\{2s_r+3p_y=s_y-12p_r=0\}\subset\mathbb P^5_{\mathbb Q}. 
\]

\begin{lemma}
\label{lem:resolved-motivic-comparison}
Given projective crepant resolutions over $\mathbb Q$,
\[
 \rho_X:\widetilde X\to X^{(0)},\qquad
 \rho_Y:\widetilde Y\to Y.
\]
For $S=\widetilde X,\widetilde Y$, there is an isomorphism
of Andr\'e motives over $\mathbb Q$,
\[
 h^3(S)\simeq\bigl(h^3(W),e_{G_{27}}\bigr)=h^3(W)^{G_{27}},\qquad
 e_{G_{27}}=\frac1{27}\sum_{h\in G_{27}}h^*.
\]
\end{lemma}

\begin{proof}
Write $S_0=X^{(0)}$ or $Y$, respectively, and
$j_X=fq_M$, $j_Y=gq_M=q_Y\sigma_W$,
using the maps of \cref{thm:quotient-H3-comparison}.
Resolve over $\mathbb Q$ the graph of $M\dashrightarrow S$:
\[
 M\xleftarrow{\ b_S\ }T_S\xrightarrow{\ a_S\ }S,
 \qquad \rho_Sa_S=j_Sb_S,
\]
with $b_S$ birational. The algebraic correspondence
\[
 \beta_S=(\sigma_W)_*(b_S)_*a_S^*:
 H^3(S)\longrightarrow H^3(W)
\]
satisfies, for $u\in H^3(S_0)$,
\[
 \beta_S\rho_S^*u
 =(\sigma_W)_*j_S^*u
 =\begin{cases}
 q_Y^*(g^*)^{-1}f^*u,&S=\widetilde X,\\
 q_Y^*u,&S=\widetilde Y,
 \end{cases}
\]
by $b_{S*}b_S^*=\sigma_{W*}\sigma_W^*=1$ and
$q_M^*g^*=\sigma_W^*q_Y^*$.
By \cref{thm:crepant-resolution-attractors,lem:rational-resolved-comparison},
$\beta_S$ identifies $H^3(S)$ with $H^3(W)^{G_{27}}$. 
We conclude by the following two commutative diagrams. 

\begin{center}
\centering
\begin{tikzpicture}[
  >=Stealth,
  font=\small,
  line width=0.45pt,
  node/.style={inner sep=2pt, outer sep=0pt},
  lab/.style={fill=white, inner sep=1pt}
]

\begin{scope}[xshift=0cm]
  \node at (5,-0.2) {$\boxed{S=\widetilde X}$};

  \node[node] (TX) at (0.0,2.1) {$T_X$};
  \node[node] (M1) at (2.2,2.1) {$M$};
  \node[node] (W1) at (5.4,2.1) {$W$};

  \node[node] (Xt) at (0.0,0.7) {$\widetilde X$};
  \node[node] (R1) at (2.2,0.7) {$M/G_{27}=R$};
  \node[node] (Y1) at (5.4,0.7) {$Y=W/G_{27}$};

  \node[node] (X0) at (0.0,-0.7) {$X^{(0)}$};

  \draw[->] (TX) -- node[lab, above] {$b_X$} (M1);
  \draw[->] (TX) -- node[lab, left] {$a_X$} (Xt);

  \draw[->] (M1) -- node[lab, right] {$q_M$} (R1);
  \draw[->] (M1) -- node[lab, above] {$\sigma_W$} (W1);
  \draw[->] (W1) -- node[lab, right] {$q_Y$} (Y1);

  \draw[->] (Xt) -- node[lab, left] {$\rho_X$} (X0);
  \draw[->] (R1) -- node[lab, below] {$f$} (X0);
  \draw[->] (R1) -- node[lab, below] {$g$} (Y1);

  \draw[->, bend right=2] (M1) to node[lab, left] {$j_X$} (X0);

  \node[align=center] at (4,-1)
    {$\displaystyle j_X=fq_M,\qquad \beta_X=\sigma_{W*}b_{X*}a_X^*$};
\end{scope}

\begin{scope}[xshift=7.5cm]
  \node at (5, -0.2) {$\boxed{S=\widetilde Y}$};

  \node[node] (TY) at (0.0,2.1) {$T_Y$};
  \node[node] (M2) at (2.2,2.1) {$M$};
  \node[node] (W2) at (5.4,2.1) {$W$};

  \node[node] (Yt) at (0.0,0.7) {$\widetilde Y$};
  \node[node] (R2) at (2.2,0.7) {$M/G_{27}=R$};
  \node[node] (Y2) at (5.4,0.7) {$Y=W/G_{27}$};

  \draw[->] (TY) -- node[lab, above] {$b_Y$} (M2);
  \draw[->] (TY) -- node[lab, left] {$a_Y$} (Yt);

  \draw[->] (M2) -- node[lab, right] {$q_M$} (R2);
  \draw[->] (M2) -- node[lab, above] {$\sigma_W$} (W2);
  \draw[->] (W2) -- node[lab, right] {$q_Y$} (Y2);

  \draw[->] (R2) -- node[lab, below] {$g$} (Y2);
  \draw[->, bend right=22] (Yt) to node[lab, below] {$\rho_Y$} (Y2);

  \draw[->] (M2) to node[lab, right] {$j_Y$} (Y2);

  \node[align=center] at (3.2,-1)
    {$\displaystyle j_Y=gq_M=q_Y\sigma_W,\qquad \beta_Y=\sigma_{W*}b_{Y*}a_Y^*$};
\end{scope}

\end{tikzpicture}
\end{center}

The motivated morphism $h^3(S)\to h^3(W)^{G_{27}}$ is an isomorphism in Betti realization. Hence it is an isomorphism of motives by faithfulness and exactness of Betti realization in Andr\'e's category \cite{andre96}. 
\end{proof}

\begin{theorem}[Modularity of the resolved attractors]
\label{thm:resolved-motivic-modularity}
Under the hypotheses of \cref{lem:resolved-motivic-comparison},
for either $S=\widetilde X$ or $\widetilde Y$ there are rank-two
Andr\'e motives $\mathsf M_{S,1},\mathsf M_{S,2}$ over $\mathbb Q$ with
\[
 h^3(S)\simeq\mathsf M_{S,1}\oplus\mathsf M_{S,2},
\]
whose Betti realizations have Hodge types
$(1,0,0,1)$ and $(0,1,1,0)$, respectively.
In particular, both resolutions are rank-two attractors.

For every prime $\ell$,
\[
 H^3_{\mathrm{\acute et}}(S_{\overline{\mathbb Q}},\mathbb Q_\ell)
 \simeq V_{1,\ell}\oplus V_{2,\ell},
 \quad
 \begin{cases}
 V_{1,\ell}=e_\Gamma H^3_{\mathrm{\acute et}}
 (W_{\overline{\mathbb Q}},\mathbb Q_\ell),\\
 V_{2,\ell}=(e_{G_{27}}-e_\Gamma)
 H^3_{\mathrm{\acute et}}(W_{\overline{\mathbb Q}},\mathbb Q_\ell),
 \end{cases}
\]
where $e_\Gamma=\frac1{24}\sum_{g\in\Gamma}g^*$.
Furthermore, with $f=f_{\mathsf{54.4.a.c}}$ and
$g=f_{\mathsf{54.2.a.a}}$, the level-$54$ newforms
of weights $4$ and $2$, respectively,
\[
 V_{1,\ell}^{\mathrm{ss}}\simeq\rho_{f,\ell}^{\mathrm{ss}},
 \qquad
 V_{2,\ell}^{\mathrm{ss}}\simeq
 \bigl(\rho_{g,\ell}(-1)\bigr)^{\mathrm{ss}}.
\]
In particular, for $p>3$, $p\ne\ell$, the Frobenius polynomial
of $H^3(S)$ is
\[
 (1-a_p(f)T+p^3T^2)(1-p\,a_p(g)T+p^3T^2).
\]
\end{theorem}

\begin{proof}
By \cref{prop:rank-two-motivic-splitting}, the image of
$e_\Gamma$ lies in that of $e_{G_{27}}$ on $H^3(W)$.
Both averaging projectors are self-adjoint for the cup-product
pairing, so
\[
 e_\Gamma e_{G_{27}}=e_{G_{27}}e_\Gamma=e_\Gamma
\]
on $h^3(W)$, by faithfulness of Betti realization.
Thus the Galois-invariant projectors
$e_\Gamma$ and $e_{G_{27}}-e_\Gamma$ cut out two
orthogonal rank-two motivated summands of
$h^3(W)^{G_{27}}$. Transporting them through
\cref{lem:resolved-motivic-comparison} proves the motivic,
Hodge and $\ell$-adic decompositions.

For $\ell=2$, the modular identifications are
\cref{thm:V1-V2-modularity}. For every $p>3$ and
$\ell\ne p$, the Frobenius traces on both summands are
independent of $\ell$: averaging the twisted Frobenius
Lefschetz traces and using \cref{lem:G27-even-cohomology}
gives precisely the formulas
\eqref{eq:G27-invariant-Frobenius-trace} and
\eqref{eq:V1-Frobenius-trace}, valid for every $\ell$.
Thus they equal $a_p(f)$ and $p\,a_p(g)$ for every $\ell$.
Chebotarev density and Brauer--Nesbitt give the stated
semisimplified isomorphisms.

Finally, each $V_{i,\ell}$ is a nondegenerate
rank-two summand for the alternating cup-product
pairing, which takes values in
$\mathbb Q_\ell(-3)$. Hence
\[
 \bigwedge\nolimits^2 V_{i,\ell}
 \simeq\mathbb Q_\ell(-3),
 \qquad
 \det V_{i,\ell}=\chi_\ell^{-3},
 \qquad i=1,2.
\]
Since $f$ and $g$ have trivial nebentypus,
Deligne's construction gives
\[
 \det\rho_{f,\ell}=\chi_\ell^{-3},
 \qquad
 \det\bigl(\rho_{g,\ell}(-1)\bigr)=\chi_\ell^{-3}, 
\]
see \cite{deligne71}; this proves the asserted
Frobenius factorization. 
\end{proof}

\appendix
\section{The point counting PARI/GP script}
\label{app:pari}

The script works with the homogeneous base coordinate $(a:b)$ and the two
Hesse cubics
\[
 a(x_0^3+x_1^3+x_2^3)=b x_0x_1x_2,
 \qquad
 b(y_0^3+y_1^3+y_2^3)=-18a y_0y_1y_2.
\]

We now describe the ${G_{27}}$-average in detail.  After dividing by the two
projective scalar subgroups, every element of ${G_{27}}$ has a unique
representative
\[
(a_0,a_1,0,b_0,b_1,0),
\qquad
a_0+a_1+b_0+b_1=0
\quad\text{in }\mathbb Z_3.
\]
Put $A=a_0+a_1$.  Over a finite extension containing $\xi_3$, choose
$c$ with $c^{p-1}=\xi_3$.  Then the fixed points of
$h\circ\operatorname{Frob}_p$ on the base and on the two projective planes
are obtained from the $\mathbb F_p$-rational points by the diagonal scalings
$c^{-A}$, $(c^{-a_0},c^{-a_1},1)$, and
$(c^{-b_0},c^{-b_1},1)$, respectively.  Substitution into the two equations
of \cref{lem:good-reduction-W} gives $\Lambda_p(g)$.

The only nontrivial step for $\Gamma$ is the solution of the semilinear
fixed-point equations.  For each of the $24$ words in the generators
$\phi,\psi$, let $\mathbf A$ denote one of the induced matrices on a projective
factor, and let $\sigma$ be the $p$-power Frobenius.  We choose an extension
$k=\mathbb F_{p^d}$ for which the twisted norm
\[
 N_d(\mathbf A)=\mathbf A \mathbf A^\sigma\cdots \mathbf A^{\sigma^{d-1}}
\]
is scalar.  Since only the fixed set of primes and the fixed $24$ group
elements are used, the admissible degrees depend only on
$p\bmod3$; they are stored in the two arrays \texttt{D1} and \texttt{D2},
so no runtime search is needed.  If $N_d(\mathbf A)=s\mathbf I$, 
multiplication of $\mathbf A$ by a
scalar does not change its projective action.  The routine \texttt{unit}
therefore chooses $\lambda\in k^\times$ with
$N_{k/\mathbb F_p}(\lambda)=s^{-1}$ and replaces $\mathbf A$ by $\lambda \mathbf A$.
After this normalization, $N_d(\mathbf A)=\mathbf I$.

The Lang--Steinberg theorem now gives a matrix $\mathbf B\in\operatorname{GL}_n(k)$ satisfying
\[
 \mathbf B^\sigma=\mathbf B\mathbf A.
\]
The program finds $\mathbf B$ by linear algebra over $\mathbb F_p$.  Indeed, for a
row vector $v\in k^n$ the equation $v^\sigma=v\mathbf A$ is $\mathbb F_p$-linear.
Writing $k$ in the polynomial basis $1,u,\ldots,u^{d-1}$ turns the map
$v\mapsto v^\sigma-v\mathbf A$ into an $nd\times nd$ matrix over $\mathbb F_p$.
Its kernel has dimension $n$; a basis of this kernel, written as rows,
forms the required matrix $\mathbf B$.  Consequently, if
$u\in\mathbb F_p^n$ and $x=\mathbf B^{-1}u$, then
\[
 \mathbf A x^\sigma=\mathbf A(\mathbf B^\sigma)^{-1}u=\mathbf B^{-1}u=x.
\]
Thus the fixed points of $\mathbf A\circ\operatorname{Frob}_p$ in projective space
are precisely $\mathbf B^{-1}\mathbb P^{n-1}(\mathbb F_p)$.  This is the set
returned by \texttt{Fix}.  Applying this construction to the two
$\mathbb P^2$ factors and the base $\mathbb P^1$, and then using
\texttt{Count}, gives $\Lambda_p(g)$; averaging over the $24$ words gives
$A_\Gamma(p)$.

\begin{mdframed}[
  linewidth=0.5pt,
  roundcorner=2pt,
  innerleftmargin=7pt,
  innerrightmargin=7pt,
  innertopmargin=5pt,
  innerbottommargin=5pt,
  skipabove=8pt,
  skipbelow=8pt
]
\footnotesize
\begin{verbatim}
default(parisize,"64M");

\\ AESZ4: G27- and Gamma-averages for a*sum(x_i^3)=b*x0*x1*x2
\\ and b*sum(y_i^3)=-18*a*y0*y1*y2.
T=[5,7,11,13,17,19,23,31,37];

\\ F_p-rational points of P^1 and P^2, viewed in the working finite field k.
P1(p,k)={my(z=0*k,o=1+0*k);concat([[o,z]~],vector(p,i,[i-1+0*k,o]~))};
P2(p,k)=
{
  my(z=0*k,o=1+0*k);
  concat([[o,z,z]~],concat(vector(p,i,[i-1+0*k,o,z]~),
    vector(p^2,n,[((n-1)\p)+0*k,((n-1)%p)+0*k,o]~)))
};
\\ C tests the two cubic equations; 
\\ Count multiplies the two fiber sizes over each fixed base point.
C(a,b,v,e)=
{
  my(v0=v[1],v1=v[2],v2=v[3],s=v0^3+v1^3+v2^3,q=v0*v1*v2);
  if(e,b*s+18*a*q==0,a*s-b*q==0)
};
Count(B,X,Y)=
{
  sum(i=1,#B,
    sum(j=1,#X,C(B[i][1],B[i][2],X[j],0)) *
    sum(j=1,#Y,C(B[i][1],B[i][2],Y[j],1)))
};

\\ G27=(Z/3Z)^3: h=[a0,a1,b0,b1,A] encodes (a0,a1,0,b0,b1,0).
\\ A=a0+a1, a0+a1+b0+b1=0 (mod 3).
m3(n)=lift(Mod(n,3));
G27elt(n)={my(r=n-1,a0=r%3,a1=(r\3)%3,b0=(r\9)%3);
  [a0,a1,b0,m3(-a0-a1-b0),m3(a0+a1)]};
G27=vector(27,n,G27elt(n));
G27avg(p)=
{
  my(d=if(p%3==1,3,2),k=ffgen(p^d,'u),xi3,c,s,B,P);
  xi3=polrootsmod(z^2+z+1,k)[1]; c=polrootsmod(x^(p-1)-xi3,k)[1];
  s=[1+0*k,c^-1,c^-2]; B=P1(p,k); P=P2(p,k);
  sum(i=1,27,my(h=G27[i],BB,X,Y);
    BB=vector(#B,j,[s[h[5]+1]*B[j][1],B[j][2]]~);
    X=vector(#P,j,[s[h[1]+1]*P[j][1],s[h[2]+1]*P[j][2],P[j][3]]~);
    Y=vector(#P,j,[s[h[3]+1]*P[j][1],s[h[4]+1]*P[j][2],P[j][3]]~);
    Count(BB,X,Y))/27
};

\\ Gamma=<phi,psi>: W lists its 24 words; 
\\ D2,D1 are the required Lang degrees for p=2,1 mod 3.
Words()=
{
  [[],[1],[2],[1,1],[1,2],[2,1],[2,2],[1,1,2],[1,2,1],[1,2,2],
   [2,1,1],[2,1,2],[2,2,2],[1,1,2,1],[1,1,2,2],[1,2,1,1],[1,2,1,2],
   [1,2,2,2],[2,1,2,1],[2,1,2,2],[1,1,2,1,1],[1,1,2,1,2],[1,2,1,2,2],[2,1,2,1,1]]
};
W=Words();
D2=[2,2,2,2,8,8,2,8,2,2,8,2,2,8,2,8,8,8,8,8,2,8,2,8];
D1=[2,3,4,3,3,3,2,6,6,6,6,6,4,4,6,4,3,6,3,6,3,4,3,4];
GammaGens(xi3)=
{
  [[xi3,0,0;0,1,0;0,0,1],[1,1,1;1,xi3,xi3^2;1,xi3^2,xi3],
   [xi3^2,0,0;0,1,0;0,0,1],[1,1,1;1,xi3^2,xi3;1,xi3,xi3^2],
   [xi3,0;0,1],[-3,1;18,3]]
};
GammaElt(v,gens)=
{
  my(E=[matid(3),matid(3),matid(2)]);
  for(i=1,#v,for(j=1,3,E[j]=E[j]*gens[2*(j-1)+v[i]])); E
};

\\ Lang--Steinberg: normalize the projective matrix, then solve v^(p)=vA over F_p.
F(M,p)=matrix(matsize(M)[1],matsize(M)[2],i,j,M[i,j]^p);
Ntw(A,p,d)={my(N=matid(matsize(A)[1]),B=A);for(i=1,d,N=N*B;B=F(B,p));N};
\\ Ntw is the twisted norm; 
\\ unit rescales A to twisted norm one without changing its PGL action.
unit(A,p,d,k)=
{
  my(B=(1+0*k)*A,s,g=ffprimroot(k),h=g^((p^d-1)/(p-1)));
  s=Ntw(B,p,d)[1,1];
  g^fflog((1+0*k)/s,h,p-1)*B
};
\\ coords uses the basis 1,k,...,k^(d-1); the kernel below is computed over F_p.
coords(a,d)={my(v=Vecrev(lift(a.pol)));concat(v,vector(d-#v,i,0))};
Lang(A,p,d,k)=
{
  my(n=matsize(A)[1],Q=vector(d,i,k^(i-1)),M=matrix(n*d,n*d),K,v,u,q);
  A=unit(A,p,d,k);
  for(c=1,n*d,
    v=vector(n,i,0*k); q=(c-1)\d+1; v[q]=Q[(c-1)%d+1];
    u=vector(n,i,v[i]^p)-v*A;
    for(i=1,n,my(t=coords(u[i],d));for(j=1,d,M[(i-1)*d+j,c]=t[j]%p))
  );
  K=matker(Mod(M,p));
  matrix(n,n,i,j,sum(r=1,d,lift(K[(j-1)*d+r,i])*Q[r]))
};

\\ If B^(p)=BA, then B^{-1}P^(n-1)(F_p) is Fix(A o Frob_p).
Fix(A,p,d,k,P)={my(B=Lang(A,p,d,k)^-1);vector(#P,i,B*P[i])};
GammaFixed(p,v,d)=
{
  my(k=ffgen(p^d,'v),xi3,E,B,P);
  xi3=polrootsmod(z^2+z+1,k)[1]; E=GammaElt(v,GammaGens(xi3)); B=P1(p,k); P=P2(p,k);
  Count(Fix(E[3],p,d,k,B),Fix(E[1],p,d,k,P),Fix(E[2],p,d,k,P))
};
GammaAvg(p)={my(D=if(p%3==2,D2,D1));sum(i=1,24,GammaFixed(p,W[i],D[i]))/24};

\\ Lefschetz gives trH3G27 and trV1; the complementary trace is trV2=trH3G27-trV1.
data(p)=
{
  my(AG27=G27avg(p),AGamma=GammaAvg(p),
     trH3G27=1+7*p+7*p^2+p^3-AG27,
     trV1=1+5*p+5*p^2+p^3-AGamma,trV2=trH3G27-trV1);
  [p,AG27,AGamma,trH3G27,trV1,trV2,trV2/p]
};
{
  printf("[%3s, %6s, %7s, %8s, %5s, %5s, %7s]\n",
  "p", "A_G27", "A_Gamma", "trH3G27", "trV1", "trV2", "trV2/p");
for(i=1,#T,
  my(r=data(T[i]));
  printf("[%3d, %6d, %7d, %8d, %5d, %5d, %7d]\n",
    r[1],r[2],r[3],r[4],r[5],r[6],r[7])
);
}
\end{verbatim}
\end{mdframed}

\vspace{1cm}

The output of the program is the following. 

\begin{mdframed}[
  linewidth=0.5pt,
  roundcorner=2pt,
  innerleftmargin=7pt,
  innerrightmargin=7pt,
  innertopmargin=5pt,
  innerbottommargin=5pt,
  skipabove=8pt,
  skipbelow=8pt
]
\footnotesize
\begin{verbatim}
[  p,  A_G27, A_Gamma,  trH3G27,  trV1,  trV2,  trV2/p]
[  5,    318,     273,       18,     3,    15,       3]
[  7,    714,     595,       22,    29,    -7,      -1]
[ 11,   2346,    2049,      -90,   -57,   -33,      -3]
[ 13,   3504,    3088,      -32,    20,   -52,      -4]
[ 17,   7128,    6516,      -72,   -72,     0,       0]
[ 19,   9588,    8866,      -68,  -106,    38,       2]
[ 23,  15996,   14754,       36,   174,  -138,      -6]
[ 31,  36534,   34705,      202,    47,   155,       5]
[ 37,  60420,   57682,       76,     2,    74,       2]
\end{verbatim}
\end{mdframed}

\newpage

\bibliography{blx}
\bibliographystyle{halpha-abbrv}

\end{document}